\documentclass[11pt,centertags]{article}

\usepackage{amsmath}
\usepackage{amsthm}
\usepackage{amscd}
\usepackage{amssymb}
\usepackage{verbatim}

\title{The Degree Theorem in higher rank}
\author{Christopher Connell~\thanks{Supported 
in part by an NSF postdoctoral fellowship.}  and Benson Farb~\thanks{Supported 
in part by the NSF and the Sloan foundation.} 
}
\date{November 7, 2000}
\newtheorem{theorem}{Theorem}[section]

\refstepcounter{theoremA}
\newtheorem*{theorem*}{Theorem}
\newtheorem{proposition}[theorem]{Proposition}
\newtheorem{lemma}[theorem]{Lemma}
\newtheorem{corollary}[theorem]{Corollary}
\newtheorem*{corollary*}{Corollary}

\newtheorem{example}[theorem]{Example}

\def\eps{\epsilon}

\def\endproof{$\diamond$ \bigskip}
\def\bar{\overline}

\newcommand\D{\partial}
\newcommand{\inner}[1]{\left\langle #1 \right\rangle }

\newcommand\op{\operatorname}

\newcommand\Jac{\operatorname{Jac}}

\renewcommand\deg{{\rm deg}}
\renewcommand\tilde{\widetilde}
\newcommand\vol{\operatorname{Vol}}
\newcommand\minvol{\operatorname{Minvol}}

\newcommand\directsum{\oplus}

\newcommand\R{\mbox{\bf R}}
\newcommand\hyp{\mbox{\bf H}}

\newcommand\tr{\operatorname{Tr}\,}

\DeclareMathOperator{\Vol}{Vol}
\DeclareMathOperator{\Isom}{Isom}
\DeclareMathOperator{\ent}{ent}
\DeclareMathOperator{\Id}{Id}
\DeclareMathOperator{\meas}{{\cal M}}
\DeclareMathOperator{\bary}{bar}
\DeclareMathOperator{\supp}{supp}
\DeclareMathOperator{\SO}{SO}
\DeclareMathOperator{\rank}{rank}
\DeclareMathOperator{\SL}{SL}

\begin{document}
\maketitle

\section{Introduction}

The problem of relating volume to degree for maps between Riemannian
manifolds is a fundamental one.  
Gromov's Volume Comparison Theorem \cite{Gr} gives such a relation for
maps into negatively curved manifolds.  In this paper we extend Gromov's
theorem to locally symmetric manifolds of nonpositive curvature.  

The key fact we prove is:  
any continuous map $f:N\rightarrow M$ between closed manifolds,
with $M$ nonpositively curved and locally symmetric (barring a few
exceptions), is homotopic to a $C^1$ map with universally bounded
Jacobian, depending only on smallest Ricci curvatures of $M,N$.  
We use this to prove the following.

\begin{theorem}[The Degree Theorem]
\label{theorem:degree}
Let $M$ be a closed, locally symmetric $n$-manifold with nonpositive
sectional curvatures.  Assume that $M$ has no local direct factors
locally isometric to $\R^k, \hyp^2$, or $\SL_3(\R)/\SO_3(\R)$.  Then for
any closed Riemannian manifold $N$ and any continuous map $f:N\to M$,
$$\deg(f)\leq C\frac{\vol(N)}{\vol(M)}$$ where $C>0$ depends only on $n$
and the smallest Ricci curvatures of $N$ and $M$.
\end{theorem}

\noindent
{\bf Remarks. }
\label{rem:theorem}
\begin{enumerate}
\item As tori have self-maps of arbitrary degree, it is easy to see
  that the theorem would be false without the ``no $\R^k$ factors''
  hypothesis. We believe that the ``no $\hyp^2$ or $\SL_3(\R)/\SO_3(\R)$
  local factors'' hypothesis is unnecessary; we show in Example \ref{ex:sl3},
  however, that the issue is delicate, and depends on the chosen
  class of metrics on $N$.

\item By scaling the metrics it is easy to see that the dependence of
  the constant $C$ on the smallest curvatures cannot be
  improved.  Actually, we determine the constant explicitly in terms
  of the symmetric space and the volume entropy of $N$ (see \S
  \ref{sec:basic}).
  
\item In \S\ref{sec:noncompact-case} we extend Theorem
  \ref{theorem:degree} to the case where $N$ and $M$ have finite
  volume (with ``bounded geometry'') but are not necessarily compact, and
  where $f$ is a coarse Lipshitz map.  \bigskip
\end{enumerate}

When $\dim(M)=2$ the conclusion of the 
theorem follows easily from the Gauss-Bonnet
Theorem.  When $\rank(M)=1$, Besson-Courtois-Gallot \cite{BCG1} proved
the stronger {\em entropy rigidity theorem}, giving the exact best
constant $C$.  Entropy rigidity is still
open in higher rank \footnote{Entropy rigidity has 
recently been proved 
\cite{BCG3,CF} for manifolds locally modelled on products of rank one
symmetric spaces with no $\hyp^2$ factors.}; 
this would correspond to the above theorem with
the constant $C$ in the inequality being
$C=\left(\frac{h(g)}{h(g_{0})}\right)^n$, where $h(g)$ is the volume
entropy of $g$ (see \cite{BCG1}), with equality being obtained iff
$g=g_{0}$ locally symmetric.  

The Besson-Courtois-Gallot technique is a central ingredient here; indeed the
main idea in our proof of Theorem \ref{theorem:degree} is to give a
higher rank version of the ``canonical map'' of \cite{BCG1}, and to give 
an {\em a priori} bound on its Jacobian.  This bound is of independent
interest, and constitutes a first step towards proving 
higher rank entropy rigidity.

\bigskip
\noindent
{\bf The Minvol invariant.} One of the basic invariants associated to a
smooth manifold $M$ is its {\em minimal volume}:
$$\minvol(M):=\inf_g\{\Vol(M,g):|K(g)|\leq 1\}$$ where $g$ ranges over
all smooth metrics on $M$ and $K(g)$ denotes the sectional curvature of
$g$.  The basic questions about $\minvol(M)$ are: for which $M$ is
$\minvol(M)>0$? when is $\minvol(M)$ realized by some metric $g$?

When a nonpositively curved manifold $M$ has a local direct factor
locally isometric to $\R^k$, 
it is easy to see that $\minvol(M)=0$.  By taking $f$
to be the identity map (while varying the metric $g$ on $M$), Theorem 
\ref{theorem:degree} immediately gives:

\begin{corollary}[Positivity of Minvol]
\label{corollary:minvol}
Let $M$ be any finite volume, 
locally symmetric $n$-manifold ($n\geq 2$) of nonpositive
curvature.  Assume that $M$ has no local direct factors 
locally isometric to $\R^k, \hyp^2$, or $\SL_3(\R)/\SO_3(\R)$.  
Then $\minvol(M)>0.$
\end{corollary}

For compact $M$, Corollary \ref{corollary:minvol} was proved
(without the $\hyp^2$ and $\SL(3,\R)$ restriction) in \cite{Gr2}
(see also \cite{Sa} for the case manifolds locally modelled on the
symmetric space for $\SL(n,\R)$).  When $M$ is real hyperbolic,
Besson-Courtois-Gallot \cite{BCG1} proved that $\minvol(M)$ is {\em
uniquely} realized by the locally symmetric metric. It seems
possible that this might hold in general.

\bigskip
\noindent {\bf Self maps and the co-Hopf property. }
As $\deg(f^n)=\deg(f)^n$, an immediate
corollary of Theorem \ref{theorem:degree} is the following.

\begin{corollary}[Self maps]
\label{corollary:hopf}
Let $M$ be a finite-volume locally symmetric manifold as in Theorem
\ref{theorem:degree}.  Then $M$ admits no self-maps of degree $>1$.
In particular, $\pi_1(M)$ is {\em co-Hopfian}: every injective
endomorphism of $\pi_1(M)$ is surjective.
\end{corollary}

Note that Corollary \ref{corollary:hopf} may also be deduced from Margulis's
Superrigidity theorem (for higher rank $M$).  The co-Hopf property for 
lattices was first proved by Prasad \cite{Pr}.

More generally, if $N$ and $M$ are as in Theorem \ref{theorem:degree}
and $f:N\to M$ and $g:M\to N$ are two maps of nonzero degree then they
both must have degree one since $f\circ g$ is a self map of $M$.

\bigskip
\noindent
{\bf Outline of the proof of Theorem
\ref{theorem:degree}. } Given 
$f:N\rightarrow M$ as in the hypothesis of the theorem, we use the
method of \cite{BCG1,BCG2} to construct a ``canonical'' map
$F:N\rightarrow M$ which is homotopic to $f$ (hence $\deg F=\deg f$).
The main difficulty is to then give a universal bound on the Jacobian of
$F$; once this is done a simple degree argument gives the theorem.

\bigskip
\noindent
{\bf Step 1 (Constructing the map):} First consider the case when the
metric on $N$ is nonpositively curved.  Denote by $Y$ (resp. $X$) the
universal cover of $N$ (resp. $M$).  Let $\meas(\partial Y),
\meas(\partial X)$ denote the spaces of atomless probability measures on
the visual boundaries of the universal covers $Y, X$.

Morally what we do is, following the method of
\cite{BCG2}, to define a map 
$$\widetilde{F}:Y\to \meas(\partial
Y)\stackrel{\phi_\ast}{\to}\meas(\partial
X)\stackrel{\bary}{\to}X$$
where $\phi_\ast=\partial
\widetilde{f}_\ast$ is the pushforward of measures and $\bary$ is the
``barycenter of a measure'' (see \S\ref{section:barycenter}). The
inclusion $Y\to \meas(\partial Y)$, denoted
$x\mapsto \mu_x$, is given by the construction of the {\em
Patterson-Sullivan measures} $\{\mu_x\}_{x\in X}$ corresponding to
$\pi_1(N)<\Isom(Y)$ (see
\S\ref{section:ps}).  An essential feature of these constructions is
that they are all canonical, so that all of the maps are {\em
  equivariant}.  Hence $\widetilde{F}$ descends to a map $F:N\to M$.

One problem with this construction outline is that the metric on $Y$ may
not be nonpositively curved.  So we must find an alternative to using
the ``visual boundary'' of $Y$.  This is done by constructing a certain
family of smooth measures $\mu_s$ on $Y$ itself, pushing them forward
via $\tilde{f}$, and convolving with Patterson-Sullivan measure on $X$.
Maps $\tilde{F}_s$ are then defined by taking the barycenters of these
measures; it is actually these maps which are considered instead of
$F$. This idea was first introduced in \cite{BCG1}.

Two new features of $F$ appear in higher rank.
First, the non-strictness of convexity of the
Busemann function (see \S\ref{section:barycenter}) must be overcome to
define $F$.  Second, and more importantly, a theorem of Albuquerque
shows that the support of each of the
measures $\mu_x$ is codimension $\rank(X)-1$ 
subset of $\partial X$ called the {\em
Furstenburg boundary} of $X$ (see \S\ref{section:ps}).  This fact and
its implications are crucial for later steps.

\bigskip
\noindent
{\bf Step 2 (The Jacobian estimate): } The heart of the paper (
\S\ref{sec:Jac} and \S\ref{section:matching}) is 
obtaining a universal bound on $F$, independent of $f$.  For this we
first realize the Jacobian of $F$ as the ratio of determinants of two
matrix integrals. We then show that whenever there are small
eigenvalues in the denominator there are a sufficient number of small
eigenvalues in the numerator with which to cancel them. The key is to find
these eigenvalues independently of the integrating measure (which
depends on $\mu_s$), therefore reducing it to a problem about
semisimple Lie groups.

\bigskip
\noindent
{\bf Step 3 (Finishing the proof): }Once a universal bound on
$|\Jac(F)|$ is found, a simple degree argument finishes the proof.  In
the case when $M$ and $N$ are not compact, the main difficulty is
proving that $F_s$ is proper.  This is quite technical, and requires
extending some of the ideas of \cite{BCS} to the higher rank setting.

\section{Patterson-Sullivan measures on symmetric spaces}
\label{section:ps}

In this section we briefly recall Albuquerque's theory \cite{Al} of 
Patterson-Sullivan measures in higher rank symmetric spaces.  
For background on nonpositively curved manifolds, symmetric spaces,
visual boundaries, Busemann functions, etc., we refer the reader to
\cite{BGS} and \cite{Eb}.  

\subsection{Basic properties}
\label{sec:basic}
Let $X$ be a Riemannian symmetric space of noncompact type.  Denote by
$\D X$ the visual boundary of $X$; that is, the set of
equivalence classes of geodesic rays in $X$, endowed with the cone
topology.  Hence $X\cup \D X$ is a compactification of $X$ which
is homeomorphic to a closed ball.  

The {\em volume entropy} $h(g)$ of a closed Riemannian $n$-manifold $(M,g)$
is defined as 
$$h(g)=\lim_{R\rightarrow \infty}\frac{1}{R} \log (\Vol(B(x,R)))$$
where $B(x,R)$ is the ball of radius $R$ around a fixed point $x$ in
the universal cover $X$.  The number $h(g)$ is independent of the
choice of $x$, and equals the topological entropy of the geodesic flow
on $(M,g)$ when the curvature $K(g)$ satisfies $K(g)\leq 0$.  
Note that while the volume $\Vol(M,g)$ is not invariant
under scaling the metric $g$, the {\em normalized entropy}
$$\ent(g)=h(g)^n\Vol(M,g)$$ is scale
invariant.

Let $\Gamma$ be a lattice in
$\Isom(X)$, so that $h(g_0)<\infty$ where $(M,g_0)$ is
$\Gamma\backslash X$ with the locally symmetric metric.

Generalizing the construction of Patterson-Sullivan, Albuquerque
constructs in \cite{Al} a family of {\em Patterson-Sullivan 
measures} on $\D X$.  This is a family of 
probability measures $\{\nu_x\}_{x\in X}$ on $\D X$ which provide a
particularly natural embedding of $X$ into the space of
probability measures on $\D X$.

\begin{theorem}[Existence Theorem, \cite{Al}]
\label{theorem:properties}
There exists a family $\{\nu_x\}$ of probability measures on $\partial
X$ satisfying the following properties:
\begin{enumerate} 

\item Each $\nu_x$ has no atoms.

\item The family of 
measures $\{\nu_x\}$ is $\Gamma$-equivariant:
$$\gamma_*\nu_x=\nu_{\gamma x} \mbox{\ for all\ }\gamma\in\Gamma$$

\item For all $x,y\in X$, the measure $\nu_y$ is absolutely
continuous with respect to $\nu_x$.  In fact the Radon-Nikodym derivative is
given explicitly by:

\begin{align}
\label{eqn:R-N}
\frac{d\nu_x}{d\nu_y}( \xi) = e^{h(g)B(x,y, \xi)}
\end{align}

\noindent
where $B(x,y,\xi)$ is the {\em Busemann function} on $X$.  For points
$x,y\in X$ and $\xi\in\D X$, the function $B:X\times X\times\D X
\rightarrow \R$ is defined by
$$B(x,y,\xi)=\lim_{t\rightarrow \infty}d_X(y,\gamma_\xi(t))-t$$
where
$\gamma_\xi$ is the unique geodesic ray with $\gamma(0)=x$ and
$\gamma(\infty)=\xi$.

\end{enumerate}
\end{theorem}

The second property implies no two measures are the same as measures.
Thus the assignment $x\mapsto \nu_x$ defines an injective map $$\nu:
X\to \meas(\D X)$$
where $\meas(\D X)$ is the space of probability
measures on $X$.  Such a mapping satisfying the above properties is
called an $h(g_0)$-{\em conformal density}.
 
\subsection{Symmetric spaces of noncompact type}

Before we present Albuquerque's theorem we will need some necessary
background about higher rank symmetric spaces.

By definition, the symmetric space $X$ is $G/K$ where $G$ is a
semisimple Lie group and $K$ a maximal compact subgroup. Fix once and
for all a basepoint $p\in X$. This choice uniquely determines a Cartan
decomposition $\mathfrak{g}=\mathfrak{k}\oplus\mathfrak{p}$ of the Lie
algebra of $G$ where $\mathfrak{k}$ is the Lie algebra of the isotropy
subgroup $K=\op{Stab}_G(p)$ of $p$ in $G$ and $\mathfrak{p}$ is
orthogonal to $\mathfrak{k}$ with respect to the killing form
$B(\cdot,\cdot)$ on $\mathfrak{g}$.  Therefore, $\mathfrak{p}$ is also
identified with the tangent space $T_pX$.

Let $\mathfrak{a}$ be, once and for all, a fixed maximal abelian
subalgebra of $\mathfrak{g}$. It follows from the Cartan decomposition
that $\mathfrak{a}\subset\mathfrak{p}$. The set
$\exp(\mathfrak{a})\cdot p$ will be a maximal flat (totally
geodesically embedded Euclidean space of maximal dimension) in $X$.
Recall, a vector $v\in TX$ is called a {\em regular vector} if it is
tangent to a unique maximal flat. Otherwise it is a {\em singular
  vector}. A geodesic is called regular (resp.  singular) if one (and
hence all) of its tangent vectors are regular (singular).  A point
$\xi\in \D X$ is regular (singular) if any (and hence all) of the
geodesics in the corresponding equivalence class are regular
(singular).

Let $\mathfrak{a}^*$ be the dual to $\mathfrak{a}$, then for each
$\alpha\in \mathfrak{a}^*$ define
$$\mathfrak{g}_\alpha=\{Y\in \mathfrak{g} | \op{ad}_A Y=\alpha(A)Y
\text{for all } A\in\mathfrak{a}\}.$$
We call $\alpha$ a {\em root} if
$\mathfrak{g}_\alpha\neq 0$. Therefore the roots form a finite set
$\Lambda$.

If $\theta_p$
is the Cartan involution associated to the point $p$, which is
$\Id$ on $\mathfrak{k}$ and $-\Id$ on $\mathfrak{p}$, then we may
define a positive definite inner product $\phi_p$ on $\mathfrak{g}$ by 
$\phi_p(Y,Z)=-B(\theta_p Y,Z)$. With respect to $\phi_p$, the
folowing {\em root space decomposition}
$$\mathfrak{g}=\mathfrak{g}_0+\sum_{\alpha\in\Lambda}\mathfrak{g}_\alpha$$
is orthogonal.

The following is proposition can be found in 2.7.3 of \cite{Eb}.

\begin{proposition}
\label{prop:roots}
Some properties of the roots and root space decomposition are:
\begin{enumerate}

\item $[\mathfrak{g}_\alpha,\mathfrak{g}_\beta]\subset
  \mathfrak{g}_{\alpha+\beta}$ if $\alpha+\beta\in\Lambda$ or is $0$
  otherwise.
\item If $\alpha\in \Lambda$ then $-\alpha\in\Lambda$ and
  $\theta_p:\mathfrak{g}_\alpha\to\mathfrak{g}_{-\alpha}$ is an
  isomorphism.
\item If $\alpha$ is not an integer multiple of some other
  $\lambda\in\Lambda$ then the only possible multiples of $\alpha$ in
  $\Lambda$ are $\pm \alpha$ and $\pm 2\alpha$.
\item We have $\mathfrak{g}_0=(\mathfrak{g}_0\cap
  \mathfrak{k})+\mathfrak{a}$.
\item If $\alpha,\beta\in\Lambda$ then
  $\beta-2\frac{\inner{\alpha,\beta}}{\inner{\alpha,\alpha}}\alpha\in
  \Lambda$ where $\inner{\cdot,\cdot}$ is the dual inner product to
  $\phi_p$ on $\mathfrak{a}^*$. Furthermore,
  $2\frac{\inner{\alpha,\beta}}{\inner{\alpha,\alpha}}$ is always an
  integer and if $\alpha$ and $\beta$ are not collinear then it is
  $\pm 1$.
\end{enumerate}
\end{proposition}

\noindent
We call a subset $\Delta\subset \Lambda$ a {\em base} for $\Lambda$ if
\begin{enumerate}

\item the elements of $\Delta$ form a basis (over $\R$) for
  $\mathfrak{a}^*$, 
\item and every root in $\Lambda$ can be written as a linear
  combination of of elements in $\Delta$ with coefficients being
  either all nonnegative integers or all nonpositive integers.
\end{enumerate}
                                
If we choose an regular element $A\in\mathfrak{a}$ then define the set
of {\em positive roots} corresponding to $A$,
$$\Lambda^+_A=\{\alpha\in\Lambda| \alpha(A)\geq 0\}$$
The subset
$\Delta^+_A\subset \Lambda^+_A$ consisting of elements which cannot be
written as a sum of two elements in $\Lambda_A^+$ is a base for
$\Lambda$. Sometimes $\Delta^+_A$ is called a {\em fundamental system
  of positive roots}.

For $A\in \mathfrak{a}$ the associated (open) {\em Weyl chamber}
$W(A)$ is the connected component of the set of regular vectors in
$\mathfrak{a}$ which contains $A$.  We also call the set $\exp
W(A)\subset \exp(\mathfrak{a})$, as well as $\exp( W(A))\cdot p\subset
X$, a Weyl chamber which we again denote by $W(A)$ using the context
to determine where exactly it lies.

The union of all the singular geodesics in the flat
$\exp(\mathfrak{a})\cdot p$ passing through $p$ is a finite set of
hyperplanes forming the boundaries of the Weyl chambers. This provides
another description of the Weyl chamber $W(A)$ as
$$W(A)=\{Y\in \mathfrak{a} | \alpha(Y)>0 \text{ for all }
\alpha\in\Delta^+_A \}.$$
For each subset $I\subset \Delta^+_A$ the
set $W_I(A)=\cap_{\alpha\in I} (\ker \alpha\cap \overline{W(A)}$ is
called the {\em Weyl chamber face} corresponding to the set $I$, and
we designate $W_{\emptyset}(A)=W(A)$. The subgroup of $K$ which
stabilizes the face $W_I(A)$ we denote by $K_I$.

\subsection{The Furstenberg Boundary}

The {\em Furstenberg boundary} of a symmetric space $X$ of noncompact
type is abstractly defined to be $G/P$ where $P$ is a minimal
parabolic subgroup of the connected component $G$ of the identity in 
$\Isom(X)$.

The Furstenberg boundary can be identified with the orbit of $G$
acting on any regular point $v(\infty)\in \D X$, the endpoint of a
geodesic tangent to a regular vector $v$. of a Weyl chamber in a fixed
flat $\mathfrak{a}$.  This follows from the fact that the action of
any such $P$ on $\D X$ fixes some regular point.

Because of this, for symmetric spaces of higher rank, behaviour on the
visual boundary can often be aptly described by its restriction to the
Furstenberg boundary.  Here we will use only some very basic
properties of this boundary.  For more details on semisimple Lie
groups and the Furstenberg boundary, see \cite{Zi}.

For a fixed regular vector $A\in\mathfrak{a}$ and associated set of
positive roots $\Lambda^+_A$ the {\em barycenter} $b$ of the Weyl
chamber $W(A)$ is defined to be
$$b=\sum_{\alpha\in\Lambda^+_A}m_\alpha H_\alpha$$
where
$m_\alpha=\dim \mathfrak{g}_\alpha$ is the multiplicty of $\alpha$ and
$H_\alpha$ is the dual vector (with respect to $\phi_p$) of $\alpha$.
Set $b^+=b/\|b\|$.

Define the set $\D_F X\subset \D X$ to be $\D_F X=G\cdot b^+(\infty)$.
Henceforth we will refer to the Furstenberg boundary as this specific
realization. We point out that for any lattice $\Gamma$ in $\Isom(X)$,
the induced action on the boundary is transitive only on $\D_F X$.
That is, $\overline{\Gamma\cdot b^+(\infty)}=G\cdot b^+(\infty)$ even
though for any interior point $x\in X$, $\overline{\Gamma\cdot x}=\D
X$.

\subsection{Albuquerque's Theorem}

Theorem 7.4 and Proposition 7.5 of \cite{Al} combine to give the
following theorem, which will play a crucial role in our proof
of Theorem \ref{theorem:degree}.
 
\begin{theorem}[Description of $\nu_x$]
\label{theorem:patterson-sullivan}
Let $(X,g_0)$ be a symmetric space of noncompact type, and let
$\Gamma$ be a lattice in $\Isom(X)$. Then
\begin{enumerate}
\item $h(g_0)=\|b\|$,
\item $b^+(\infty)$ is a regular point, and hence $\D_F X$ is a
  regular set,
\item For any $x\in X$, the support $\supp(\nu_x)$ of $\nu_x$ is equal
  to $\D_F X$, and
\item $\nu_x$ is the unique probability measure invariant under the
  action on $\D_F X$ of the compact isotropy group $\op{Stab}_G(x)$ at
  $x$. In particular, $\nu_p$ is the unique $K$-invariant probability
  measure on $\D_F X$.
\end{enumerate}
\end{theorem}

Note that when $X$ is has rank one, $\D_FX=\D X$. In general
$\D_F X$ has codimension $\rank(X)-1$ in $\D X$.

\subsection{Limits of Patterson-Sullivan measures}
We now describe the asymptotic behaviour of the $\nu_x$ as $x$ tends
to a point in $\D X$.

For any point $\xi$ of the visual boundary, let $S_\theta$ be the set
of points $\xi\in \D_F X$ such that there is a Weyl chamber $W$ whose
closure $\D \bar{W}$ in $\D X$ contains both $\theta$ and $\xi$. Let
$K_\theta$ be the subgroup of $K$ which stabilizes $S_\theta$.
$K_\theta$ acts transitively on $S_\theta$ (see the proof below).

\begin{theorem}[Support of $\nu_x$]
\label{theorem:support}
Given any sequence $\{x_i\}$ tending to $\theta\in \D X$ in the cone
topology, the measures $\nu_{x_i}$ converge in $\meas(\D_F X)$ to the
unique $K_\theta$-invariant probability measure $\nu_{\theta}$
supported on $S_\theta$.
\end{theorem}

\begin{proof} 
  Let $x_i=g_i \cdot p$, for an appropriate sequence $g_i\in G$.
  Recall that $\nu_{x_i}=(g_i)_* \nu_p$. Then combining part (4) of
  Theorem~\ref{theorem:patterson-sullivan} with Proposition 9.43 of
  \cite{GJT} have that some subsequence of the $\nu_{x_i}$ converges
  to a $K_\theta$-invariant measure $\nu_\theta$ supported on
  $S_\theta$. 
  
  Note that in \cite{GJT}, the notation $I$ refers to a subset of a
  fundamental set of roots corresponding to the face of a Weyl chamber
  containing $\theta$ in its boundary. If $g_i\cdot p=k_i a_i\cdot p$
  converges then both $k=\lim k_i$ and $a^I=\lim_i a_i^I$ exist (note
  the definition of $a^I$ in \cite{GJT}). Again in the notation of
  \cite{GJT}, $K_\theta$ is the conjugate subgroup $(ka^I) K^I
  (ka^I)^{-1}$ in $K$. Moreover, $S_\theta$ is the orbit $k a^I
  K^I\cdot b^+(\infty)$.
  
  By Corollary 9.46 and Proposition 9.45 of \cite{GJT} any other
  convergent subsequence of the $\nu_{x_i}$ produces the same measure
  in the limit, and therefore the sequence $\nu_{x_i}$ itself
  converges to $\nu_\theta$ uniquely.
\end{proof}

In the case when $\theta$ is a regular point, the above theorem
implies that $S_\theta$ is a single point and the limit measure
$\nu_\theta$ is simply the Dirac probability measure at that point
point in $\D_FX$.

\section{The barycenter of a measure}
\label{section:barycenter}

In this section we describe the natural map which is an essential 
ingredient in the method of Besson-Courtois-Gallot.  

Let $\phi$ denote the lift to universal covers of $f$ with
basepoint $p\in Y$ (resp. $f(p)\in X$), i.e. $\phi=\tilde{f}:Y\to X$.
We will also denote the metric and Riemannian volume form on
universal cover $Y$ by $g$ and $dg$ respectively.  Then for each
$s>h(g)$ and $y\in Y$ consider the probability measure $\mu_y^s$ 
on $Y$ in the Lebesgue
class with density given by
$$\frac{d\mu_y^s}{dg}(z)=\frac{e^{-sd(y,z)}}{\int_{Y}
  e^{-sd(y,z)}dg}.$$
The $\mu_y^s$ are well defined by the choice of $s$.

Consider the push-forward $\phi_*\mu_y^s$, which is a measure on $X$.
Define $\sigma_y^s$ to be the convolution of $\phi_*\mu_y^s$ 
with the Patterson-Sullivan measure $\nu_z$ for the symmetric metric. 

In other words, for $U\subset\partial X$ a Borel set, define 
$$\sigma_y^s(U)=\int_{X}\nu_z(U)d(\phi_*\mu_y^s)(z)$$

Since $\Vert \nu_z \Vert=1$, we have
$$\Vert\sigma_y^s\Vert=\Vert\mu_y^s\Vert=1.$$

Let
$B(x,\theta)=B(\widetilde{f}(p),x,\theta)$ be the Busemann function on
$X$ with respect to the basepoint $\widetilde{f}(p)$ (which we
will also denote by $p$).  For $s>h(g)$ and $x\in X, y\in Y$ define a 
function
$$\mathcal{B}_{s,y}(x)=\int_{\D X}B(x,\theta)d\sigma_y^s(\theta)$$

By Theorem \ref{theorem:support}, the support of $\nu_z$, hence of
$\sigma_y^s$, is all of $\partial_FX$, which in turn equals the
$G$-orbit $G\cdot b^+(\infty)$. Hence
$$\mathcal{B}_{s,y}(x)=\int_{\D_F
  X}B(x,\theta)d\sigma_y^s(\theta)=\int_{G\cdot
  b^+(\infty)}B(x,\theta)d\sigma_y^s(\theta).$$

Since $X$ is nonpositively curved, the Busemann function $B$ is
(non-strictly) convex on $X$.  Hence $\mathcal{B}_{s,y}$ is convex on
$X$, being a convex integral of convex functions.  While $B$ is strictly
convex only when $X$ is negatively curved, we have the following.

\begin{proposition}[Strict convexity of $\mathcal{B}$]
  For each fixed $y$ and $s$, the function $x\mapsto
  \mathcal{B}_{y,s}(x)$ is strictly convex, and has a unique critical
  point in $X$ which is its minimum.
\end{proposition}

\begin{proof}
  It suffices to show that given a geodesic segment $\gamma(t)$
  between two points $\gamma(0),\gamma(1)\in X$, there exists some
  $\xi\in \D_F X$ such that function $B(\gamma(t),\xi)$ is {\em
    strictly} convex in $t$, and hence on an open positive
  $\mu_y$-measure set around $\xi$. We know it is convex by the
  comment preceding the statement of the proposition.
  
  If $B(\gamma(t),\xi)$ is constant on some geodesic subsegment of
  $\gamma$ for some $\xi$, then $\gamma$ must lie in some flat
  $\mathcal{F}$ such that the geodesic between $\xi\in\partial
  \mathcal{F}$ and $\gamma$ (which meets $\gamma$ at a right angle) also
  lies in $\mathcal{F}$. On the other hand, $\xi\in \partial_F X$ is in
  the direction of the algebraic centroid in a Weyl chamber, and
  $\gamma$ is perpendicular to this direction. By the properties of the
  roots, $\gamma$ is a regular geodesic (i.e. $\gamma$ is not contained
  in the boundary of a Weyl chamber).  In particular, $\gamma$ is
  contained in a unique flat $\mathcal{F}$.  Furthermore, $\partial_F
  X\cap \D \mathcal{F}$ is a finite set (an orbit of the Weyl group). As
  a result, for almost every $\xi\in \D_F X$ $B(\gamma(t),\xi)$ is
  {\em strictly} convex in $t$.

  For fixed $z\in X$, by the last property listed in Theorem 
  \ref{theorem:properties}, we see that
  $$\int_{\D_F X} B(x,\theta)d\nu_z(\theta)$$
  tends to $\infty$ as
  $x$ tends to any boundary point $\xi\in\D X$. Then for fixed $y$ and
  $s>h(g)$, $\mathcal{B}_{y,s}(x)$ increases to $\infty$ as $x$ tends
  to any boundary point $\xi\in\D X$. Hence it has a local minimum in
  $X$, which by strict convexity must be unique.  
\end{proof}

We call the unique critical point of $\mathcal{B}_{s,y}$ the {\em
barycenter} of the measure $\sigma_y^s$, and define a map 
$\widetilde{F}_s:Y\to X$ by
$$\widetilde{F}_s(y)= \text{the unique critical point of
$\mathcal{B}_{s,y}$}$$ 
Since for any two points $p_1,p_2\in X$
$$B(p_1,x,\theta)=B(p_2,x,\theta)+B(p_1,p_2,\theta)$$
we see that
$\mathcal{B}_{s,y}$ only changes by an additive constant when we
change the basepoint of $B$. Also, $\mathcal{B}_{s,y}$ only changes by
a multiplicative constant when we change the basepoint in the
definition of $\mu_y$. Since neither change affects the critical point
of $\mathcal{B}_{y,s}$, we see that $\widetilde{F}_s$ is independent
of choice of basepoints.

The equivariance of $\widetilde{f}$ and of $\{\mu_y\}$ implies that
$\widetilde{F}_s$ is also equivariant. Hence $\widetilde{F}_s$
descends to a map $F_s: N\to M$. It is easy to see that $F_s$ is
homotopic to $f$.

\begin{proposition}
\label{prop:homotopy}
The map $\Psi_s\colon [0,1]\times N\to M$ defined by
$$\Psi_s(t,y)=F_{s+\frac t{1-t}}(y)$$
is a homotopy between
$\Psi_s(0,\cdot)=F_s$ and $\Psi_s(1,\cdot)=f$.
\end{proposition}

\begin{proof}
  From its definitions, $\tilde{F}_s(y)$ is continuous in $s$ and $y$.
  Observe that for fixed $y$, $\lim_{s\to
    \infty}\sigma^s_y=\nu_{\phi(y)}$. If follows that $\lim_{s\to
    \infty}\tilde{F}_s(y)=\phi(y)$. This implies the proposition.
\end{proof}

As in \cite{BCG1}, we will see that $F_s$ is $C^1$, and will estimate
its Jacobian.

\section{The Jacobian estimate}
\label{sec:Jac}

Let $X$ be expressed as a product of its irreducible factors 
$X=X_1\times \ldots\times X_k$, and let $g_i$ denote the restricted symmetric
metric on each factor $X_i$. As above, $h(g_i)$ denotes the volume entropy
of $(X_i,g_i)$. The main estimate of this paper is the following.

\begin{theorem}[The Jacobian Estimate]
\label{JacF}
For all $s>h(g)$ and all $y\in N$ we have
$$\vert \Jac F_s(y)\vert\le C\left(\frac{s}{h(g_{1}) h(g_2)\cdots
    h(g_k)}\right)^n$$
for some constant $C$, depending only on $\dim
M$.
\end{theorem}

\bigskip
\noindent
{\bf Dependence of constants. } Up to scaling of the metric, there are
only a finite number of irreducible symmetric spaces of noncompact
type in a given dimension. Therefore it is sufficient to show that $C$
depends only on the individual symmetric spaces $(X_i,g_{i})$.
Furthermore, when we apply Theorem \ref{JacF}, we will take the limit
as $s\to h(g)$ so that the quantity
$C\left(\frac{h(g)}{h(g_{1})h(g_2)\cdots h(g_k)}\right)^n$ is the
constant appearing in Theorem \ref{theorem:degree}. It is evident then 
that the right hand side of inequality of Theorem \ref{theorem:degree} 
is scale invariant with respect to the metrics $g$ and $g_i$.

We claim that the quantities $h(g)$ and $h(g_{i})$ can be bounded by
Ricci curvatures. The Bishop volume comparison Theorem (\cite{BC})
states that if the Ricci curvatures of $(Y,g)$ are all greater than
$(n-1)\kappa$ for some $\kappa\leq 0$ then for any $y\in Y$ and $r>0$,
$$\vol B(y,r)\leq V_{\kappa}(r)$$
where $V_{\kappa}(r)$ is the volume
of the ball of radius $r$ in the space form of constant curvature
$\kappa$. In particular this implies that $$h(g)\leq
\lim_{r\to\infty}\frac{\log V_{\kappa}(r)}{r}=(n-1)\sqrt{-\kappa}.$$

Similarly, in the course of the proof of Theorem \ref{JacF} we will
see explicitly that $$h(g_{i})=\tr \sqrt{-R_i(b^+,\cdot,b^+,\cdot)}$$
where $R_i$ is the curvature tensor on $(X_i,g_i)$. In particular
$$h(g_i)\geq \min\{1,-\op{Ricci}(b^+,b^+)\}.$$
Therefore the constant
$C$ in Theorem \ref{theorem:degree} depends only on the Ricci
curvatures of $N$ and $M$.

\bigskip

We will prove Theorem \ref{JacF} in several steps.  

\subsection{Finding the Jacobian}

We obtain the differential of $F_s$ by implicit differentiation:
$$0=D_{x=F_s(y)}\mathcal{B}_{s,y}(x)=\int_{\D_F X}
d{B}_{(F_s(y),\theta)}(\cdot) d\sigma^s_y(\theta)$$

Hence as $2$-forms
\begin{gather*}
0=D_y D_{x=F_s(y)}\mathcal{B}_{s,y}(x)=\int_{\D_F X}
Dd{B}_{(F_s(y),\theta)}(D_y F_s(\cdot),\cdot) d\sigma^s_y(\theta)\\
 -s \int_{Y}\int_{\D_F X} d{B}_{(F_s(y),\theta)}(\cdot)
\inner{ \nabla_y d(y,z),\cdot} d\nu_{\phi(z)}(\theta) d\mu^s_y(z)
\end{gather*}

The distance function $d(y,z)$ is Lipschitz and $C^1$ off of the cut
locus which has Lebesgue measure 0. It follows from the Implicit
Function Theorem (see \cite{BCG2}) that $F_s$ is $C^1$ for $s>h(g)$.
By the chain rule,
$$\Jac F_s = s^n\frac{\det\left(\int_{Y}\int_{\D_F X}
    d{B}_{(F_s(y),\theta)}(\cdot) \inner{ \nabla_y
      d(y,z),\cdot} d\nu_{\phi(z)}(\theta) d\mu^s_y(z)
  \right)}{\det\left(\int_{\D_F X}
    Dd{B}_{(F_s(y),\theta)}(\cdot,\cdot)
    d\sigma^s_y(\theta)\right)}$$

Applying H{\"o}lder's inequality to the numerator gives:

$$|\Jac F_s|\leq s^n\frac{\det\left(\int_{\D_F X}
    d{B}_{(F_s(y),\theta)}^2 d\sigma^s_y(\theta)\right)^{1/2}
  \det\left(\int_{Y} \inner{\nabla_y d(y,z),\cdot}^2
    d\mu^s_y(z)\right)^{1/2}}{\det\left(\int_{\D_F X}
    Dd{B}_{(F_s(y),\theta)}(\cdot,\cdot)
    d\sigma^s_y(\theta)\right)}$$

Using that $\tr \inner{\nabla_y d(y,z),\cdot}^2 =\left|\nabla_y
  d(y,z)\right|^2=1$, except possibly on a measure $0$ set, we may
estimate
$$\det\left(\int_{Y} \inner{\nabla_y d(y,z),\cdot}^2
  d\mu^s_y(z)\right)^{1/2}\leq \left(\frac{1}{\sqrt{n}}\right)^n$$

Therefore
\begin{align}
\label{eq:Jac1}
| \Jac F_s|\leq
\left(\frac{s}{\sqrt{n}}\right)^n\frac{\det\left(\int_{\D_F X}
    d{B}_{(F_s(y),\theta)}^2
    d\sigma^s_y(\theta)\right)^{1/2}}{\det\left(\int_{\D_F X}
    Dd{B}_{(F_s(y),\theta)}(\cdot,\cdot)
    d\sigma^s_y(\theta)\right)}
\end{align}

\subsection{Reduction to Irreducible Case}

In this subsection we make, following \cite{CF}, a reduction to 
the case when $X=\widetilde{M}$ is irreducible.  

If $X=X_1\times \ldots\times X_k$ is the irreducible expression for
$X$ as a product, the group $G=\Isom(X)$ can also be written as a product
$G=G_1\times G_2\cdots \times G_{k}$, where each $G_i\neq 
\op{SL}(2,\R),\op{SL}(3,\R)$ is a simple Lie group.  Theorem
~\ref{theorem:patterson-sullivan} implies that for all $y\in Y$, the
measure $\sigma^s_y$ is supported on the $G$-orbit $$G\cdot
b^+(\infty)=\{(G_1\times G_2\cdots \times G_{k})\cdot b^+(\infty)\}$$

Hence 
$$\D_FX=G\cdot b^+(\infty)=\D_FX_1\times\cdot\times \D_FX_{k}$$
Since each $X_i$ has rank one, $\D_FX_i=\D X_i$ so that 
$$\D_FX=\D X_1\times\cdots\times \D X_{k}$$

Let $B_i$ denote the Busemann function for the rank one symmetric
space $X_i$ with metric $g_i$. Then for $\theta_i\in \D X_i\subset \D
X$ and $x,y\in X_i$ we have
$B(x,y,\theta_i)=B_i(x,y,\theta_i)$. Since the factors
$X_i$ are orthogonal in $X$ with respect to the metric $g_{0}$, the
Busemann function of $(X,g_{0})$ with basepoint $p\in X$ at a point
$\theta=(\theta_1,\dots,\theta_{k})\in\D_F X$ is given by
$$B(x,\theta)=\frac{1}{\sqrt{{k}}}\sum_{i=1}^{k} B_i(x_i,\theta_i).$$

The Schur estimate for the determinant of symmetric semidefinite block
matrices states,
$$\det \begin{pmatrix} A & B \cr B^* & C \end{pmatrix}\le
\det(A)\det(C).$$
Applying the dual form of this estimate to our
symmetric tensors we have
\begin{gather*}
  \det\left(\int_{\D_FX} \left(\sum_{i=1}^{k}
      d(B_i)_{(\pi_i F_s(y),\pi_i\theta)}\right)^2
    d\sigma^s_y(\theta)\right)\leq\\
  \prod_{i=1}^{k} \det \left(\int_{\D_FX_i}
    \left(d(B_i)_{(\pi_i F_s(y),\theta_i)}\right)^2
    d(\pi_{i})_{*}\sigma^s_y(\theta_i)\right),
\end{gather*}
where $\pi_i:X\to X_i$ and $\pi_i:\D_F X\to \D_F X_i$ are the
canonical projections.

Since $Dd{B}_{(F_s(y),\theta)}=\frac{1}{\sqrt{{k}}}\sum_{i=1}^{k}
Dd{B_i}_{(\pi_i F_s(y),\pi_i\theta)},$ the denominator already splits as,
$$\det\left(\int_{\D_F X} Dd{B}_{(F_s(y),\theta)}(\cdot,\cdot)
  d\sigma^s_y(\theta)\right) = \prod_{i=1}^{k} \det \left(\int_{\D_FX_i}
  \left(Dd(B_i)_{(\pi_i F_s(y),\theta_i)}\right)
  d(\pi_{i})_{*}\sigma^s_y(\theta_i)\right).$$

Putting these together we obtain,
\begin{gather*}
  |\Jac F_s(y)|\leq \left(\frac{s}{\sqrt{n}}\right)^{n}
  \prod_{i=1}^{k} \frac{ \det \left(\int_{\D_FX_i}
      \left(d(B_i)_{(\pi_i F_s(y),\theta_i)}\right)^2
      d(\pi_{i})_{*}\sigma^s_y(\theta_i)\right)^{1/2}}{ \det
    \left(\int_{\D_FX_i} \left(Dd(B_i)_{(\pi_i
          F_s(y),\theta_i)}\right)
      d(\pi_{i})_{*}\sigma^s_y(\theta_i)\right)}.
\end{gather*}
Therefore we only need to bound each term in the product seperately.
It suffices then to prove that for an irreducible symmetric space
$(X,g_0)\neq \hyp^2, \op{SL}(3,\R)/\op{SO}(3,\R)$, 
and for any measure $\mu$ on $\D_FX$, that 
$$\frac{\det\left(\int_{\D_F X}
    d{B}_{(F_s(y),\theta)}^2
    d\mu(\theta)\right)^{1/2}}{\det\left(\int_{\D_F X}
    Dd{B}_{(F_s(y),\theta)}(\cdot,\cdot)
    d\mu(\theta)\right)}\leq \frac{C}{h(g_0)}.$$

We will continue to write $\sigma_y^s$ instead of $\mu$ or
$(\pi_{i})_{*}\sigma^s_y$. The only property we use of $\sigma_y^s$
from this point on is that it is fully supported on $\D_F X$. Since
$\op{supp}((\pi_{i})_{*}\sigma^s_y)=\pi_i(\op{supp}(\sigma^s_y))=\D_F
X_i$ there is no harm by this imprecision.

\subsection{Simplifying the Jacobian}
\label{sec:simpJac}

As stated above we need only now consider irreducible $(X,g_{0})$.  For each
point $x\in X$, we let $\mathcal{F}_x$ denote the canonical flat
passing through $x$, i.e. $\mathcal{F}_x=\exp(\mathfrak{a})\cdot x$.
We denote the tangent space to $\mathcal{F}_x$ simply as $\mathcal{F}$
with the base point suppressed since it is naturally isomorphic to the
Lie algebra $\exp(\mathfrak{a})$.

We wish to bound the quantity
$$\frac{\det\left(\int_{\D_F X}
    d{B}_{(F_s(y),\theta)}^2
    d\sigma^s_y(\theta)\right)^{1/2}}{\det\left(\int_{\D_F X}
    Dd{B}_{(F_s(y),\theta)}(\cdot,\cdot)
    d\sigma^s_y(\theta)\right)}$$

Let $\mathcal{F}$ denote the tangent space to the flat
$\mathcal{F}_{F_s(y)}$.  Choose an orthonormal basis $\{ e_i\}$ for
the tangent space $T_{F_s(y)}X$ such that $e_1,\dots,e_{\rank(X)}$ is
a basis for $\mathcal{F}$ with $e_1(\infty)=b^+(\infty)$.  We may
write the term

\begin{equation}
\label{eq:denom}
\int_{\D_F X}Dd{B}_{(F_s(y),\theta)}(\cdot,\cdot)
d\sigma^s_y(\theta)
\end{equation}

in matrix form as $$\int_{\D_F X} O_{\theta}
\begin{pmatrix} 0 & 0 \\ 0 & D_\lambda \end{pmatrix}O_\theta^* 
\ d\sigma^s_y(\theta)$$
where $O_\theta$ is the orthogonal matrix in
the $e_i$ basis corresponding to the derivative of the unique isometry
in $K=\op{Stab}_G(F_s(y))$ which sends $e_1$ to $v_{(F_s(y),\theta)}$
(the vector in the tangent space of the point $F_s(y)$ in the
direction $\theta\in\D_F X$). In the above expression, the upper left
zero matrix sub-block has dimensions $\rank(X)\times \rank(X)$, and
$D_\lambda$ has the form
$$D_\lambda=\begin{pmatrix} \lambda_1 & 0 &
\cdots & 0 \\ 0 & \lambda_2 & 0 & \vdots \\ \vdots & 0 & \ddots & 0 \\ 0
& \cdots & 0 & \lambda_{n-\rank(X)}
\end{pmatrix}$$
where $\{\lambda_1,\ldots ,\lambda_{n-\rank(X)}\}$ is the set of
nonzero eigenvalues of $Dd{B}_{(F_s(y),\theta)}$. Since
$Dd{B}_{(x,\theta)}$ is $G$ equivariant, its eigenvalues do not depend
on $x$ but only on which $K$-orbit in $\D X$ the point $\theta$ lies
in. In particular, $Dd{B}_{(x,\theta)}$ is flow invariant and hence
the Ricatti equation shows that it is simply related to the curvature
tensor by
$$Dd{B}_{(x,\theta)}=\sqrt{-R(v_{(x,\theta)},\cdot,v_{(x,\theta)},\cdot)}$$

On the other hand in a symmetric space
$R(v,\cdot,v,\cdot)=-(\op{ad}_v)^2|_{\mathfrak{p}}$. Therefore the
eigenvalues of $Dd{B}_{(F_s(y),\theta)}$ are those of
$Dd{B}_{(p,b^+(\infty))}$ which in turn are those of
$\sqrt{\op{ad}_{b^+}^2}|_{\mathfrak{p}}$. (Note that while
$\op{ad}_{b^+}$ does not preserve $\mathfrak{p}$,
$(\op{ad}_{b^+})^2|_{\mathfrak{p}}$ is a symmetric endomorphism of
$\mathfrak{p}$.) Recall, $b^+=b/\|b\|$ where
$b=\sum_{\beta\in\Lambda^+_A} m_\beta H_\beta$ for any choice of
$A\in\mathfrak{a}$ (the choice of $A$ only determines the Weyl chamber
containing $b$). Setting 
$$\mathfrak{p}_\alpha=\mathfrak{p}\cap
(\mathfrak{g}_\alpha\directsum\mathfrak{g}_{-\alpha}),$$
we have $\mathfrak{p}_\alpha=\{X-\theta_p X: X\in
\mathfrak{g}_{\alpha} \}$.

By definition of $\mathfrak{g}_{\alpha}$, for each
$\alpha\in\Lambda^+_A$ we may write
$$(\op{ad}_{b^+})^2|_{\mathfrak{g}_{\alpha}}=\alpha(b^+)^2\Id=
\left(\frac{1}{\|b\|}\sum_{\beta\in\Lambda^+_A} \alpha(
  m_\beta H_\beta)\right)^2\Id.$$
The same expression clearly holds for
$(\op{ad}_{b^+})^2|_{\mathfrak{g}_{-\alpha}}$. Therefore, for any
$\alpha\in \Lambda$,
$\sqrt{(\op{ad}_{b^+})^2|_{\mathfrak{p}_{\alpha}}}=|\alpha(b^+)|$. For
$\mathfrak{p}_0=\mathfrak{a}$ the same formula holds with $\alpha=0$.
In particular, the ratio of the largest eigenvalue (denoted by
$\lambda_{\max}$) among the $\lambda_i$'s in $D_\lambda$ to the
smallest nonzero eigenvalue (denoted by $\lambda_{\min}$) only depends
on $X$.

Furthermore, since $\alpha(b^+)>0$ for all $\alpha\in \Lambda^+_A$ and
$\dim \mathfrak{p}_\alpha=m_\alpha$, we have
\begin{gather*}
  \tr\sqrt{\op{ad}_{b^+}^2|_{\mathfrak{p}}}
  =\sum_{\alpha\in\Lambda^+_A} m_\alpha\alpha(b^+)=\frac{1}{\| b\|}
  \sum_{\alpha,\beta\in\Lambda^+_A} m_\alpha m_\beta\alpha(
  H_\beta)=\phantom{\hspace{90pt}}\\
  \phantom{\hspace{120pt}}\frac{1}{\| b\|}\inner{
    \sum_{\beta\in\Lambda^+_A}m_\beta H_\beta, \sum_{\alpha\in\Lambda^+_A}
    m_\alpha H_\alpha}=\frac{\| b\|^2}{\| b\|}=h(g_0)
\end{gather*}
where the last equality follows from Theorem
\ref{theorem:patterson-sullivan}. As a result, there is a constant $c$
only depending on $X$ such that
\begin{equation}
\label{eq:evals}
\frac{h(g_{0})}{c}\leq \lambda_i \leq c\  h(g_{0})
\end{equation}
for $i=1,\ldots ,(n-\rank(X))$.  We now use the following.

\begin{lemma}
\label{lemma:detsum}
  The determinant of a sum of $n\times n$ positive semidefinite
  matrices is a nondecreasing homogeneous polynomial of degree $n$ in
  the eigenvalues of each summand. Furthermore, if the sum is positive 
  definite, then the determinant is strictly increasing in the
  eigenvalues of the summands.
\end{lemma}

\begin{proof}
  Let $M$ be the sum of positive semidefinite matrices. Then there
  exist fixed orthogonal matrices $O_l$ and real numbers $\lambda_{l,j}$ 
such that $M$ may be written as 
  $$M=\sum_l O_l\begin{pmatrix} \lambda_{l,1} & 0 & \cdots & 0 \\ 0 &
    \lambda_{l,2} & 0 & \vdots \\ \vdots & 0 & \ddots & 0 \\ 0 &
    \cdots & 0 & \lambda_{l,n}
  \end{pmatrix} O_l^*$$ 
  Then we have the differentiation formula (see, e.g. Prop. 2.8 of
\cite{Ch}): 
  $$\frac{d}{d\lambda_{l,j}}\det M=\tr \left(\frac{d}{d\lambda_{l,j}}
    M\right)M^{\op{adj}}$$
  where $M^{\op{adj}}$ is the adjunct matrix
  of $M$. Now,
  $$\frac{d}{d\lambda_{l,j}} M = O_l E_{(j,j)} O_l^*$$
  where
  $E_{(j,j)}$ is the elementary matrix with $1$ in the $(j,j)$
  position and zeros elsewhere. Therefore, by cyclically permuting
  $O_l$ in the trace above we find that $\frac{d}{d\lambda_{l,j}}\det
  M$ is the $(j,j)$ the entry of $O_l^* M^{\op{adj}} O_l$ which is
  nonnegative since $M$ is positive semidefinite. Lastly, if $M$ is
  positive definite then $O_l^* M^{\op{adj}} O_l$ is also, which means
  that $\frac{d}{d\lambda_{l,j}}\det M$ is positive. The lemma
  follows.
\end{proof}

Applying Lemma \ref{lemma:detsum} to the Riemann sums for the 
integral (\ref{eq:denom}) above, using the bound in Equation
(\ref{eq:evals}), and taking limits, gives 

\begin{gather*}
\begin{split}
\det\ \int_{\D_F X} &Dd{B}_{(F_s(y),\theta)}(\cdot,\cdot)
  d\sigma^s_y(\theta)\geq \\
&\left(\frac{h(g_{0})}{c}\right)^n\det\ \int_{\D_F X}
  O_{\theta}\begin{pmatrix} 0 & 0 \\ 0 & I_{n-\rank(X)}
  \end{pmatrix}O_\theta^* \ d\sigma^s_y(\theta)
\end{split}
\end{gather*}
where $I_{n-\rank(X)}$ is the identity matrix of dimension
$n-\rank(X)$.

Next we observe that, relative to the orthonormal basis $\{e_1,\ldots
,e_{\rank(X)}\}$ for $T_{F_s(y)}X$, the expression 
$$\int_{\D_F X} d{B}_{(F_s(y),\theta)}^2
  d\sigma^s_y(\theta)$$
may be written in the form
$$Q_1=\int_{\D_F X} O_{\theta}\begin{pmatrix} 1 & 0_{1\times (n-1)} \\ 
  0_{(n-1)\times 1} & 0_{(n-1)\times(n-1)} \end{pmatrix}O_\theta^*
\ d\sigma^s_y(\theta)$$
where $O_\theta$ is the same matrix as above.  Let 
$$Q_2= \int_{\D_F X}
  O_{\theta}\begin{pmatrix} 0 & 0 \\ 0 & I_{n-\rank(X)}
  \end{pmatrix}O_\theta^* \ d\sigma^s_y(\theta)$$

We have just shown that, to prove Theorem \ref{JacF}, it suffices to
prove that 
\begin{equation}
\label{eq:mainbound}
\frac{\displaystyle \det Q_1}{\displaystyle (\det Q_2)^2}\leq C
\end{equation}
for some constant $C$.  The rest of this section will be devoted to
proving this.

\subsection{Eigenvalue matching}

Here is the general idea of our proof of Theorem \ref{JacF}, which we
have reduced to showing (\ref{eq:mainbound}) above.  
Since the numerator is bounded
above, we consider when the matrix $Q_2$ in the denominator has any
eigenvalues smaller than a certain constant depending only on the
dimension of $X$. When this occurs, Theorem \ref{prop:det-est} below
will show that each such eigenvalue is matched by at least two smaller
(up to a universal constant) eigenvalues of the matrix $Q_1$ in the
numerator. 

Let $\{v_i\}$ be an orthonormal eigenbasis for the symmetric matrix
$Q_2$, and recall that $\{e_i\}$ is a basis for the tangent space $\cal
F$ to the 
fixed, chosen flat.  Note that the $i$-th eigenvalue of the matrix $Q_2$ may be
written as $$L_i=v_i^{\ast}Q_2v_i=\int_{\D_F X}\sum_{j=\rank(X)+1}^n
\inner{O_{\theta}.e_j,v_i}^2\ d\sigma^s_y(\theta).$$

We first argue that no $L_i$ equals zero. Since $s>h(g)$ we have that
the measures $\mu_y^s$ is a finite measure in the Lebesgue class
($dg$). Since the $\nu_x$ for $x\in X$ are positive on any open set
(with respect to the cone topology) of $\D_F X$, it follows that
$\sigma_y^s$ is as well. In particular, $\{O_\theta| \theta\in
\op{supp}(\sigma_y^s)=\D_F X\}$ is isomorphic to the group $K$ and
therefore there is no nonzero subspace $V\subset T_{F_s(y)}X$ such
that $O_\theta V\subset \mathcal{F}$ for all $\theta\in\D_F X$. Hence
none of the eigenvalues $L_i$ are $0$.

Let $\epsilon=1/(\rank(X)+1)$.  Note that $\epsilon$ is a constant
depending only on $n$, as there are only finitely many symmetric
spaces of a given rank and given dimension.  Suppose $k$ of the
eigenvalues are strictly less than $\epsilon$.  Since each $L_i\leq
1$, and since $$\sum_i L_i=\tr Q_2=n-\rank(X)$$
it follows easily that
$k\leq \rank(X)$.  By rearranging the order we may assume that
$L_i<\eps$ for $i=1,\cdots,k$.

Let $H$ be an inner product space over $\R$, and denote by $\SO(H)$ the
special orthogonal group of $H$. Scale the bi-invariant metric on 
$\SO(H)$ so that $\SO(H)$ has diameter $\pi/2$.  Define the {\em angle} 
between two subspaces $V,W\subset
H$ as $$\angle (V,W):=\inf\{d_{\SO(H)}(I,P): P\in \SO(H) \text{ with } 
PV\subset W\text{ or } PW\subset V\} $$

Let $\pi_V(W)$ represent the orthogonal projection of $W$ onto $V$.
Then it is routine to verify the following properties of the angle: 
\begin{enumerate}
\label{angleproperties}
\item $\angle (V,W)\leq \frac{\pi}2$
\item $\angle (V,W)= \angle (W^\perp,V^\perp$)
\item $\angle (V,W)=\angle (W,V)$
\item If $V\subseteq U$ and $\dim U\leq \dim W$ then $\angle (V,W)\leq
  \angle (U,W)$, or \\ if $V\subseteq U$ and $\dim V\geq \dim W$ then
  $\angle (V,W)\geq \angle (U,W)$
\item If $\angle V,W=0$ then $V\subseteq W$ or $W\subseteq V$
\item If $U\subseteq W$ then $\angle (\pi_V(U),U)\leq \angle
  (\pi_V(U),W)\leq\angle (V,W)$
\end{enumerate}

For a $1$-dimensional subspace $V$ spanned by a vector $v$, our 
definition of angle agrees with the usual definition:
\begin{enumerate}
\item[7.] $V=\text{span}\{v\}\Rightarrow 
\cos(\angle(V,W))=\frac{\inner{v,\pi_W(v)}}{|v|\cdot|\pi_v(W)|}$
\end{enumerate}

Finally, $\angle$ satisfies the following form of the triangle
inequality.

\begin{lemma}[Triangle inequality for $\angle$]
\label{lemma:triangle}
Let $U,V,W$ be subspaces of a fixed inner product space $H$.  Suppose
that $\dim U=\dim W\leq \dim V$.  Then
$$\angle(V,W)\leq \angle(U,V)+\angle(U,W)$$
\end{lemma}

\begin{proof}
By definition of $\angle$ there exist $P_1,P_2,P_3\in \SO(H)$ with 
\begin{itemize}
\item $P_1W\subseteq V$ and $\angle(V,W)=d_{\SO(H)}(I,P_1)$.
\item $P_2U\subseteq V$ and $\angle(U,V)=d_{\SO(H)}(I,P_2)$.
\item $P_3U=W$ and $\angle(U,W)=d_{\SO(H)}(I,P_3)$.
\end{itemize}

Now $P_2P_3^{-1}W\subseteq V$ so that 
$$
\begin{array}{rl}
d(I,P_1)&\leq d(I,P_2P_3^{-1})\\
&\\
&=d(P_2,P_3)\\
&\\
&\leq d(I,P_2)+d(I,P_3)
\end{array}
$$
and we are done.

\end{proof}

One of the main ingredients in the proof of Theorem \ref{JacF} is the
following.

\begin{theorem}[Eigenvalue Matching Theorem]
\label{prop:det-est}
  For any $k$-frame given by orthonormal vectors $v_1,\dots,v_k$ of
  $T_{x}X$ with $k\leq\rank(X)$ there is an orthonormal $2 k$-frame
  given by vectors $v'_1,v''_1\dots,v'_{k},v''_{k}$, each perpendicular
  to $\op{span}\{v_1,\dots,v_k\}$, such that for
  $i=1,\dots,k$ and all $h\in K$, there is a constant $C$, depending
  only on $\dim X$, such that
  $$\angle (hv'_{i},\mathcal{F}^\perp) \leq C\angle (hv_i,\mathcal{F})$$
  and
  $$\angle (hv''_{i},\mathcal{F}^\perp) \leq C\angle (hv_i,\mathcal{F})$$
  where $hv$
  represents the linear (derivative) action of $K$ on $v\in T_{x}X$.
\end{theorem}

We will prove Theorem \ref{prop:det-est} in the next section; its
proof is independent of the rest of the paper.

\subsection{Proof of  the Jacobian Estimate}

Assuming Theorem \ref{prop:det-est} for the moment, we now complete
the proof of Theorem ~\ref{JacF}.

\bigskip
\noindent
{\bf Proof of Theorem \ref{JacF}. }
  From equation \eqref{eq:Jac1} and the reduction in \S
  \ref{sec:simpJac} we see that it is sufficient to show that
  $$\frac{\det Q_1}{(\det Q_2)^2}\leq C$$
  for some constant $C$
  depending only on $n$.
  
  As before let $L_1,\ldots,L_k$ be the $k\leq \rank(X)$ eigenvalues
  of $Q_2$ which are strictly less than $\epsilon=1/(\rank(X)+1)$.  If
  no such eigenvalues exist, then there is a lower bound on $Q_2$
  depending only on $\rank(X)$.  As there is an upper bound on $Q_1$,
  we are done (see the discussion on dependency of constants above).
  So we assume $k\geq 1$.
  
  Let $v_1,\ldots,v_k$ be an orthonormal set of associated
  eigenvectors.  Recall that $\{ e_i\}$ denotes the chosen orthonormal
  basis for the $T_{F_s(y)}X$ such that $e_1,\dots,e_{\rank(X)}$ spans
  the tangent space $\mathcal{F}$ to the fixed maximal flat.
  
  For any vector $v\in T_{F_s(y)}X$ let
  $$r(v)=\sum_{j=\rank(X)+1}^n \inner{e_j,v}^2$$
  so that
  $$L_i=\int_{\D_F X} r(O_\theta^*v_i)\ d\sigma^s_y(\theta).$$
  
  Since $e_1,\ldots,e_{\rank(X)}$ form an orthonormal basis for
  $\mathcal{F}$, for any unit vector $v$ we have
  $$
  \begin{array}{ll}
  \cos (\angle
  (v,\mathcal{F}))&=\inner{v,\pi_\mathcal{F}(v)}/|\pi_\mathcal{F}(v)|\\ 
  &=\inner{v,\sum \inner{v,e_j}e_j}/(\sum \inner{v,e_j}^2)^{1/2}\\
  &=(\sum\inner{v,e_j}^2)^{1/2}
  \end{array}
  $$
  so that
  $$\cos (\angle (v,\mathcal{F}))^2=\sum_{j=1}^{\rank(X)}\inner{v,e_j}^2$$
  
  Hence
  \begin{align*}
    r(v)&=1-\sum_{j=1}^{\rank(X)} \cos^2(\angle v,e_j)\\
    &=1- \cos^2(\angle v,\mathcal{F})\\
    &=\sin^2(\angle v,\mathcal{F})
  \end{align*}
  
  Similarly
  $$\inner{v,e_1}^2\leq
  \sum_{j=1}^{\rank(X)}\inner{v,e_j}^2=\sin^2(\angle
  v,\mathcal{F}^\perp)$$
  
  For each $i=1,\ldots,k$, let $v_{i}'$ and $v_{i}''$ be the pair of
  vectors corresponding to $v_i$ produced by the Eigenvalue Matching
  Theorem (Theorem \ref{prop:det-est}).  That theorem together with
  the concavity of $\sin^2 \theta$ for $0\leq\theta\leq \pi/2$ gives,
  for all $\theta\in \D_F X$ and for each $w_i=v'_i$ or $v''_i$, that
  $$\sin^2(\angle O_\theta^*w,\mathcal{F}^\perp)\leq \sin^2(C \angle
  O_\theta^*v_i,\mathcal{F})\leq C \sin^2(\angle O_\theta^*
  v_i,\mathcal{F})$$
  where $C>1$ is the constant in the Eigenvalue
  Matching Theorem.
  
  Furthermore, $Q_1$ is the integral (against a probability measure)
  of matrices with all eigenvalues less than $1$ so no eigenvalue of
  $Q_1$ is greater than one. Hence we may estimate

  \begin{align*} \det Q_1&\leq \prod_{i=1}^{k}
    (v'_{i}.Q_1.v'_i)(v''_i.Q_1.v''_i)\\
    &=\prod_{i=1}^{k}\left(\int_{\D_F
        X}\inner{e_1,O_{\theta}^*.v'_i}^2\ 
      d\sigma^s_y(\theta)\right)\left(\int_{\D_F
        X}\inner{e_1,O_{\theta}^*.v''_i}^2\ d\sigma^s_y(\theta)\right)\\
    &\leq\prod_{i=1}^{k}\left(\int_{\D_F X}\sin^2 (\angle
      O_{\theta}^*.v'_i,\mathcal{F}^\perp)\ 
      d\sigma^s_y(\theta)\right)\left(\int_{\D_F X}\sin^2 (\angle
      O_{\theta}^*.v''_i,\mathcal{F}^\perp)\ 
      d\sigma^s_y(\theta)\right)\\ &\leq
    \prod_{i=1}^{k}\left(\int_{\D_F X} C\sin^2 (\angle
      O_{\theta}^*.v_i,\mathcal{F})\ 
      d\sigma^s_y(\theta)\right)\left(\int_{\D_F X}C\sin^2 (\angle
      O_{\theta}^*.v_i,\mathcal{F})\ d\sigma^s_y(\theta)\right)\\ 
    &=C^k \prod_{i=1}^{k} L_i^2\\ &= C^k \det Q_2^2 \prod_{i=k+1}^{n}
    L_i^{-2}\\ &\leq C^k \det Q_2^2 (\rank(X)+1)^{2(n-k)}
  \end{align*}
  
  The last inequality follows from the definition of $k$, whereby
  $L_i\geq\frac{1}{\rank(X)+1}$ for each $i>k$.
  
  The constant $C$ in Theorem
  \ref{JacF} may be taken to be the product (over factors $X_j$ of $X$ with
  dimension $n_j$),
  $$\frac{1}{\sqrt{n}^n} \prod_{j} C_j^{\rank(X)/2}
  c_j^{n_j}(\rank(X_j)+1)^{(n_j)}$$
  where $C_j\geq 1$ is the constant
  $C$ from Theorem \ref{prop:det-est}, $c_j$ is the constant $c$ in
  Equation (\ref{eq:evals}) and $k_j$ is the constant $k$ above. This
  combined constant depends only on $n=\dim X$.  \endproof

\subsection{A cautionary example}
\label{ex:sl3}
In the general method of \cite{BCG1} as well as
here, one is solving a minimization problem without regard to the
measure.  However, at least in the $\SL_3(\R)/\SO_3(\R)$ case, 
to get a bound on the Jacobian of $F_s$ one must use further properties
of the measure, as indicated by the example we now give.

If for a single flat $\mathcal{F}_x$ and a
sequence of $y_i\in \mathcal{F}_x$, the measures $\sigma_{y_i}^s$ tend
to the sum of Dirac measures 
$\frac{1}{2}\delta_{b^+(\infty)}+\frac{1}{2}\delta_{w b^+(\infty)}$ 
where $w$ is in the Weyl group for $\mathcal{F}_x$, then we claim that
$\Jac F_s(y_i)\to\infty$. First note that the sum
$$dB_{(F_s(y_i),b^+(\infty))}^2+dB_{(F_s(y_i),w b^+(\infty))}^2$$ has
only a $3$-dimensional kernel, while
$$DdB_{(F_s(y_i),b^+(\infty))}+DdB_{(F_s(y_i),w b^+(\infty))}$$ has 
a $2$-dimensional kernel. Furthermore 
$$Q_1=\int_{\D_F X} dB_{(F_s(y_i),\theta)}^2
d\sigma_{y_i}^s\ \ \mbox{and}\ \ Q_2=\int_{\D_F X} DdB_{(F_s(y_i),\theta)}
d\sigma_{y_i}^s$$ degenerate in the same way, so that 
$\det(Q_1)/\det(Q_2)^2$ is unbounded. This can be easily verified explicitly
in the case of a sum of five Dirac measures for which both integrals
are nonsingular degenerating to the sum of the two Dirac measures
given above.

A similar problem occurs when there are $\hyp^2$ factors.

\section{Proof of the 
Eigenvalue Matching Theorem}
\label{section:matching}

In order to prove Theorem \ref{prop:det-est} we will need a series
of lemmas.

\subsection{Dimension inequalities}

For any $x\in X$ and any subspace $V\subseteq T_xX$, denote by $K_V$
the elements of $K$ which stabilize $V$ (i.e. leave $V$ invariant).
For $V\subset \mathcal{F}$, if $\op{Fix}_K(V)$ is the subgroup of $K$
which fixes $V$ pointwise then $K_V=U \cdot \op{Fix}_K(V)$ where $U$
is the subgroup stabilizing $V$ of the (discrete) Weyl group which
stabilizes $\mathcal{F}$ (see \cite{Eb}).

The following lemma is a basic algebraic ingredient in the proof 
of Theorem \ref{prop:det-est}.

\begin{lemma}[Dimension inequality, I]
\label{lem:dim}
With the above notations, $$\dim \left(\op{span} \{K_V\cdot
\mathcal{F}\}^\perp \right)\geq 2\dim(V).$$
\end{lemma}

\begin{proof}
  First we show that $K_V\cdot \mathcal{F}$ is itself a subspace hence
  equal to its span.
  
  Recognize that $K_V \cdot \mathcal{F}$ is the union of all tangent
  spaces to flats which contain $V$. Pick a basis $v_1,\ldots,v_l$ of $V$
  note that $K_V\cdot \mathcal{F}=\cap_{i=1}^l \mathcal{F}(v_i)$ where
  $\mathcal{F}(v_i)$ is the union of all the tangent spaces to flats
  containing $v_i$ using the notation of \cite{Eb}. Proposition 2.11.4
  of \cite{Eb} states that $\mathcal{F}(v_i)=\R^r\times X_i$ for some
  symmetric space of noncompact type and $r\leq \rank(X)$. In
  particular it is a manifold and the tangent space to it corresponds
  to $K_{v_i}\cdot \mathcal{F}$, which is a vector space. Then
  $K_V\cdot \mathcal{F}$ is a vector space.
  
  Let $K_{\mathcal{F}}$ be the stabilizer of $\mathcal{F}$ in $K$.
  Then $K_{\mathcal{F}}\subset W\cdot K_V$ where $W$ denotes the Weyl
  group (a finite group).  Hence $\dim
  K_{\mathcal{F}}=\dim(K_{\mathcal{F}}\cap K_V)$.  Hence
  $$\dim K_V\cdot \mathcal{F}=\dim K_V + \dim \mathcal{F} -\dim
  K_{\mathcal{F}}$$
  
  Since $X=K\cdot \mathcal{F}$ we obtain $$\dim M=\dim K +\dim
  \mathcal{F} -\dim K_{\mathcal{F}}.$$
  Putting this together we
  obtain, $$(\dim \op{span} \{K_V\cdot F\})^\perp=\dim M -\dim
  K_V\cdot \mathcal{F}=\dim K-\dim K_V.$$
  
  But Lemma \ref{lem:dimK} below gives that this final term is $\geq
  2\dim V$, as desired.
\end{proof}

The following lemma was used in the proof of Lemma \ref{lem:dim}.  
Recall that, at this point, we are assuming that the symmetric space 
$X$ is irreducible and has $\rank(X)\geq 2$.

\begin{lemma}[Dimension inequality, II]
\label{lem:dimK}
Assume that $X\neq \SL_3(\R)/\SO_3(\R)$. Then for any subspace $V\subset
\mathcal{F}$, we have
$$\dim K\geq 2 \dim V +\dim K_V.$$
\end{lemma}

This lemma is the only place where $X\neq \SL_3(\R)/\SO_3(\R)$ is used.

\begin{proof}
  For a root $\alpha\in\Lambda$ in $\mathcal{F}$, define
  $\mathfrak{k}_\alpha=(\Id+\theta_p)\mathfrak{g}_\alpha$, where
$\theta_p$ is the Cartan involution at $p=F_s(y)$. Then by
  Proposition 2.14.2 of \cite{Eb} we have that
  $\mathfrak{k}_\alpha=\mathfrak{g}_\alpha\directsum\mathfrak{g}_{-\alpha}\cap
  \mathfrak{k}$, $\mathfrak{k}_\alpha=\mathfrak{k}_{-\alpha}$, and
  $\dim \mathfrak{k}_\alpha\geq 1$.
  
  Note that from the definition of $\mathfrak{g}_\alpha$ it follows
  immediately that 
  $$\mathfrak{k}_\alpha=\{Y\in \mathfrak{k} | [X,Y]=0 \text{ for all }
  X\in \ker\alpha\}.$$
  
  Note that in $G$ the normalizer mod centralizer is finite for any
  flat subspace. Therefore for any $V\subset\mathcal{F}$ we may write
  the Lie algebra $\mathfrak{k}_V$ of $K_V$ as,
  $$\mathfrak{k}_V=\{Y\in \mathfrak{k} | [X,Y]=0 \text{ for all } X\in
  V\}.$$
  It then follows from the previous statements that,
  $$\mathfrak{k}_V=\mathfrak{k}_0+\sum_{\stackrel{\alpha\in\Lambda}
    {V\subset\ker\alpha}} \mathfrak{k}_\alpha.$$
  
 Consequently, we may assume that $V$ in the statement of the lemma is
  {\em maximally singular}: $V$ may be written as the intersection of
  the kernels of the greatest number of roots among all subspaces of
  dimension $\dim V$.  Otherwise $\dim K_V=\dim \mathfrak{k}_V$ is
  strictly smaller than it would be if $V$ were maximally singular.
  
  Recall that we have the invariant inner product $\phi_p$ on
  $\mathfrak{a}$ and hence on $\mathcal{F}$. Let $\Lambda$ denote the
collection of roots.  For $\alpha\in\Lambda$, let 
$H_\alpha\in \mathcal{F}$ denote the dual root vector (with respect
  to $\phi_p$) corresponding to $\alpha$.  For any subset $V\subset
  \mathcal{F}$ we define the function
  $$\op{card}_R(V):=\frac{1}{2}\op{card}\{\alpha\in\Lambda |
  H_\alpha\in V\}.$$
  Since root vectors lying in a subspace always
  come in opposing pairs, $\op{card}_R$ is a positive integer.
  
Let $\alpha$ be any root.  Note that if a subspace $V\subset\ker\alpha$, then
  $H_\alpha$ lies in $V^\perp$.  Therefore the statement of the lemma
  reduces to showing that
  $$\dim \mathfrak{k}_0+\sum_{\alpha\in\Lambda}
  \dim\mathfrak{k}_\alpha \geq 2\dim V +\dim
  \mathfrak{k}_0+\sum_{\stackrel{\alpha\in\Lambda}
    {V\subset\ker\alpha}} \dim \mathfrak{k}_\alpha,$$
  or more simply,
  $$\sum_{H_\alpha\in \mathfrak{F}\setminus V^\perp}\dim
  \mathfrak{k}_\alpha \geq 2\dim V.$$
  
  Swapping $V^\perp$ for $V$ and vice versa, and using $\dim
  \mathfrak{k}_\alpha\geq 1$ for each $\alpha$, 
it is sufficient to prove that
  \begin{gather}
    \label{eq:dim}
    \op{card}_R(\mathcal{F}\setminus V)\geq 2(\rank(X)-\dim V).
  \end{gather}
  
  Since we are assuming that $G$ is simple, we could check this
  condition by using a classification of root vectors in the simple
  algebras such as in \cite{Va}. However, because this would be
  tedious we will instead give a synthetic proof.
  
  For each $i=0,\ldots,\op{rank}(X)$, we say that $W_i\subset
  \mathcal{F}$ is a {\em maximally rooted subspace of dimension} $i$ if
  $$\op{card}_R(W_i)=\max \{\op{card}_R(V) : V\subset
  \mathcal{F}\text{ with }\dim V=i\}.$$
  In other words, $W_i$ is
  maximally rooted if $W_i^\perp$ is maximally singular.  We claim
  that if $0=W_0,W_1,\ldots,W_{\rank(X)}=\mathcal{F}$ are any
  maximally rooted subspaces of $\mathcal{F}$ with $\dim W_i=i$, then
  for $0<i\leq \rank(X)$,
\begin{equation}
\label{eq:cardind}
\op{card}_R(W_i)\geq i+\op{card}_R(W_{i-1})
\end{equation}
  
  This is true for $i=1$ since $W_1$ is one dimensional it contains a
  root vector pair and the trivial subspace $W_0$ contains none. By
  induction, assume the claim holds for all maximally rooted subspace
  $W_i$ of dimension $i<j$.  In particular, for such a space $W_{j-1}$
  and for any subspace $Z\subset W_{j-1}$ of codimension one,
  $\op{card}_R(Z)\leq \op{card}_R(W_{j-2})$ so
  $$\op{card}_R(W_{j-1}\setminus
  Z)=\op{card}_R(W_{j-1})-\op{card}_R(Z)\geq j-1.$$
  We claim that there exists a
  root vector $H_\alpha$ which is not in $W_{j-1}$ or its
  perpendicular $W_{j-1}^\perp$ (with respect to $\phi_p$). If
  not, then every root vector either lies in $W_{j-1}$ or $W_{j-1}^\perp$
  which implies the root system is reducible (e.g. Corollary 27.5 of
  \cite{Hu2}), and hence $G$ is reducible, contrary to
  assumption.
  
  Therefore, $H_{\alpha}^\perp\cap W_{j-1}$ is a codimension one
  subspace of $W_{j-1}$ and by inductive hypothesis there are at least
  $j-1$ distinct pairs of root vectors $\pm H_{\alpha_1},\ldots,\pm
  H_{\alpha_{j-1}}$ in $W_{j-1}\setminus (H_{\alpha}^\perp\cap
  W_{j-1})$. For each of these we have
  $\phi_p(H_\alpha,H_{\alpha_l})\neq 0$.  By the standard calculus of
  roots (e.g. Proposition 2.9.3 of \cite{Eb}) this implies that for
  each $1\leq l \leq j-1$ either $\pm (H_\alpha + H_{\alpha_l})$ or
  $\pm(H_{\alpha}-H_{\alpha_l})$ is a pair of root vectors lying in
  $W_{j-1}\directsum\inner{H_{\alpha}}$ which does not lie in
  $W_{j-1}$.  Including $H_{\alpha}$, these form at least $j$ pairs of
  root vectors which are contained in
  $W_{j-1}\directsum\inner{H_{\alpha}}\setminus W_{j-1}$. Therefore
  $\op{card}_R(W_{j-1}\directsum\inner{H_{\alpha}})\geq
  \op{card}_R(W_{j-1})+j$. Since by definition of $W_j$,
  $\op{card}_R(W_{j})\geq
  \op{card}_R(W_{j-1}\directsum\inner{H_{\alpha}})$, the claim follows.
  
  Recursively applying Equation \ref{eq:cardind} shows that for $0\leq i
  < j\leq \rank(X)$, $$\op{card}_R(W_j)-\op{card}_R(W_i)\geq
  \sum_{k=i}^{j}k=\frac{j(j+1)}{2}-\frac{i(i+1)}{2}.$$
  
  Now to prove the inequality~\eqref{eq:dim}, as noted before we may
  assume $V$ of dimension $q$ is maximally rooted, since then
  $V^\perp$ is maximally singular. Since $\mathcal{F}$ is a maximally
  rooted space, the above expression reads 
  $$\op{card}_R(\mathcal{F}\setminus
  V)=\op{card}_R(\mathcal{F})-\op{card}_R(V)=\frac{\rank(X)(\rank(X)+1)}{2}
  -\frac{q(q+1)}{2}.$$
  
  This is readily seen to be greater that $2(\rank(X)-q)$ unless
  $\rank(X)=2$ and $q=0$ ($V=\mathcal{F}$). However, every irreducible
  Lie algebras of rank two other than $\mathfrak{sl}(3,\R)$ has at
  least four pairs of roots (see \cite{Hu1}, p.44, Figure 1),
  and hence the inequality \eqref{eq:dim} is satisfied in all of the
  required cases.
\end{proof}

\subsection{Angle inequalities}

\begin{lemma}[Angle inequality, I]
\label{lemma:angle1}  
For any subspace $V\subseteq \mathcal{F}$ there is a subspace $V'\subset
V^\perp$ with $\dim V'\geq 2 \dim V$ and a constant $C$ depending only
on the symmetric space $X$ such that for all $k\in K$,
$$\angle (kV', \mathcal{F}^\perp) \leq C\angle (kV,\mathcal{F})$$
where $kV$ represents
the linear (derivative) action of $K$ on $V\subset T_{x}X$.
\end{lemma}

\begin{proof} 
  For any subspace $V\subset \mathcal{F}$, let
  $U_1,U_2,\ldots,U_{l(V)}$ be the maximally singular subspaces of
  dimension $\dim V$ which have minimal angle with $V$. Define
  $S_V=U_1\directsum\ldots\directsum U_{l(V)}\subset \mathcal{F}$. If
  $G(r,\mathcal{F})$ denotes the Grassmann variety of subspaces in
  $\mathcal{F}$ with dimension $r$, then the set of $V\in
  G(r,\mathcal{F})$ for which $l(V)$ is constant has codimension
  $l(V)-1$ in $G(r,\mathcal{F})$.
  
  For any subspace $V\subset \mathcal{F}$ we define a subspace 
$V'\subset \mathcal{F}^\perp$ by
$$V'=(\op{span}\{ K_{S_V}\cdot \mathcal{F}\})^\perp$$
where $K_{S_V}$ is the subgroup of $K$ which stabilizes
  $S_V$.  By Proposition \ref{lem:dim}, $V'$ has dimension at least
  $2\dim V$ since we always have $K_{S_V}\subset K_U$ for some
  $U\subset \mathcal{F}$ with $\dim U=\dim V$.
  
  If no such constant $C$ as in the lemma exists then there is a
  sequence $k_i\in K$ and $V_i\subset \mathcal{F}$ with $\dim V_i=r$
  such that
  \begin{gather}
  \label{eq:angle}
  \frac{\angle (k_i\,V_i,\mathcal{F})}{\angle (k_i\,V_i',
    \mathcal{F}^\perp)}\to 0.
  \end{gather}
  
  Now since $S_V$ and hence $V'$ varies upper semicontinuously in $V$
  (thinking of the map $V\rightarrow V'$ as a self-map of
  $G(r,\mathcal{F})$), it follows from the continuity of the $\angle$
  function that
  $$\frac{\angle (k\,V,\mathcal{F})}{\angle (k\,V',
    \mathcal{F}^\perp)}$$
  is lower semicontinuous in $V$.
  
  However since both $K$ and $G(r,\mathcal{F})$ are compact, for some
  subsequence of the $k_i V_i$, the $k_i$ converge to $k_0\in K$ and
  the $V_i$ converge to a fixed subspace $V_0\subset \mathcal{F}$.
  Furthermore, $k_0 V_0$ lies in $\mathcal{F}$ since $\angle
  (k_0V_0,\mathcal{F})$ must be $0$. It follows that $k_0\in W\cdot
  K_{V_0}$ where $W$ is the Weyl group stabilizing $\mathcal{F}$.
  
  By construction, $K_{V_0}\subset K_{V_0'}$ and for any $w\in W$,
  $$\angle (w V_0',\mathcal{F}^\perp)=\angle
  (V_0',w^{-1}\mathcal{F}^\perp)= \angle (V_0',\mathcal{F}^\perp).$$
  
  Therefore, we also have $\angle (k_0V_0',\mathcal{F}^\perp)=0$.
  Continuity of $\angle$ along with the fact that $W\subset K$ acts
  isometrically implies that it is sufficient to show that for any
  fixed subspace $V\subset \mathcal{F}$ the quantity $$\liminf_{k\to
    K_V}\frac{\angle (kV,\mathcal{F})}{\angle (kV',
    \mathcal{F}^\perp)}$$
  is bounded away from $0$. Note that since this quantity is lower
  semicontinuous in $V$, and since $G(r,\mathcal{F})$ is compact, it is 
  unnecessary to show that the bound is independent of $V$.

  First we handle the denominator. Using the bi-invariance of the metric
  on $\SO(n)$, the properties of the angle function, and the fact that
  for all $k_0\in K_{S_V}$ we have $k_0\,k\,k^{-1}\mathcal{F}\subset
  K_{S_V}\mathcal{F}$, it follows that 

\begin{align*}
d_{\SO(n)}(k,K_{S_V})= d_{\SO(n)}(k^{-1},K_{S_V})&=d_{\SO(n)}(K_{S_V}\cdot k,
Id)\\ &\geq
\inf\{d_{\SO(n)}(\Id,P): P\in \SO(n) \text{ with } P
k^{-1}\mathcal{F}\subset K_{S_V}\mathcal{F} \}\\ &=\angle
(\span\{K_{S_{V}}\mathcal{F}\},k^{-1}\mathcal{F}) \\ &=\angle (k K_{S_V}
\mathcal{F},\mathcal{F}) \\ &=\angle \left((k K_{S_V}
\mathcal{F})^\perp,\mathcal{F}^\perp\right) \\ &=\angle (k
V',\mathcal{F}^\perp) \\ 
\end{align*}
                      
  So it remains to show that for any sequence $k_i\to k^V\in K_V$ in
  any fixed neighborhood $U$ of $K_V$, that $\angle (k_i V,\mathcal{F})
  \geq C d_{\SO(n)}(k_i,K_{S_V})$. Furthermore, since $\angle (k_i
  V,\mathcal{F})=\angle (k_i (k^V_i)^{-1}V,\mathcal{F})$ 
for any $k^V_i\in K_V$, we 
  may assume that $k_i\to \Id$.
  
  By Theorem 2.10.1 of \cite{Va}, in a sufficiently
  small neighborhood of $\Id$ we may uniquely write $k_i$ as $k_i=\exp({\rm
    k}^\perp_i)\exp({\rm k}^S_i)$ where ${\rm
    k}^S_i\in\mathfrak{k}_{S_{V_0}}$ and ${\rm
    k}^\perp_i\in\mathfrak{k}_{S_{V_0}}^\perp$.
  Furthermore ${\rm k}^S_i\to 0$ and ${\rm k}^\perp_i\to 0$.
  
  Bi-invariance of the metric on $\SO(n)$ implies that for $|{\rm
    k}^\perp_i|<\frac{\pi}{2},$
  $$d_{\SO(n)}(k_i,K_{S_V})=d_{\SO(n)}(\exp({\rm
    k}^\perp_i),K_{S_V})=|{\rm k}^\perp_i|.$$
  
Now $K_V$ is the only subgroup of $K$ which both 
leaves $V$ in $\mathcal{F}$ and also intersects all sufficiently small
neighborhoods of the identity.  Therefore, in order
  to show that $\angle (k_i V,\mathcal{F}) \geq C |{\rm k}_i^\perp|$, we 
need only show that 
$$d_{\SO(n)}(k_i,K_V)/|k_i^\perp| \not \rightarrow 0$$

Well, the Cambell-Baker-Hausdorff formula implies that 
  $$\exp({\rm k}^\perp_i)\exp({\rm k}^S_i)=\exp\left({\rm
      k}^\perp_i+{\rm k}^S_i+O(|{\rm k}^\perp_i|\cdot |{\rm
      k}^S_i|)\right).$$
  Since the definition of $S_V$ implies that
  $\mathfrak{k}_{S_V}\supset \mathfrak{k}_V$ and ${\rm k}^\perp_i$ is
  perpendicular to $\mathfrak{k}_{S_V}$, we have
  $$d_{\SO(n)}(k_i,K_V)\geq |{\rm k}_i^\perp|+O(|{\rm k}^\perp_i|\cdot
  |{\rm k}^S_i|).$$
  Since we had $|{\rm k}^S_i|\to 0$ this finishes
  the lemma.

\end{proof}

\begin{lemma}[Angle inequality, II]
\label{lem:subspace}
For any subspace $V$ of $T_{x}X$ with $\dim V\leq \rank(X)$, there is
a subspace $V'\perp V$ with $\dim V'\geq 2 \dim V$, and a constant $C$
depending only on $n$, such that
$$\angle (kV',\mathcal{F}^\perp )\leq C\angle (kV,\mathcal{F}) \ \ \ 
\mbox{for all}\ k\in K$$
\end{lemma}

\begin{proof} 
The first step of the proof is to reduce to the case when $V$ is a
subspace of $\mathcal{F}$, so that Lemma \ref{lemma:angle1} may be applied.  

We first observe that the lemma is true if and only if it is 
true with $V$ replaced by $k_0V$ for any fixed $k_0\in K$.  Since $K$ is 
compact we may therefore choose $V$ among all $kV, k\in K$ so that 
$\angle (V,\mathcal{F})\leq \angle (kV,\mathcal{F})$ for all $k\in K$.
  
With this assumption, consider the projection $W=\pi_F(V)$ of $V$ onto
$\mathcal{F}$.  By Lemma \ref{lemma:angle1}, we obtain a subspace $W'$
such that
$$\angle (kW',\mathcal{F}^\perp)\leq C\angle (kW,\mathcal{F})$$ 
for all $k\in K$. Then we let $V'$
be the projection of $W'$ onto $V^\perp$. By the properties of the angle
function (see \ref{angleproperties}), it follows that

$$
\begin{array}{rll}
\angle (kV',\mathcal{F}^\perp)&\leq \angle
(kW',\mathcal{F}^\perp)+\angle (kV',kW')& 
\mbox{by Lemma \ref{lemma:triangle}}\\
&&\\
&\leq C\angle (kW',\mathcal{F}^\perp)+ \angle (V',W')&\\
&&\\
&\leq C\angle (kW,\mathcal{F}) +\angle (V',W')&\mbox{since\
  }(W^{\perp})^\perp\supseteq W\\
&&\\
&\leq C\angle (kW,\mathcal{F}) +\angle V,\mathcal{F}&\mbox{for same reason}\\
&&\\
&=C\angle (kW,\mathcal{F})+\angle (V,\mathcal{F})&\mbox{since\ }W=\pi_F(V)
\end{array}
$$

Thus it suffices to bound $\angle (kW,\mathcal{F})$ by a constant times
$\angle(kV,\mathcal{F})$.  But 

$$
\begin{array}{rll}
\angle(kW,\mathcal{F})&\leq
\angle(kV,\mathcal{F})+\angle(kV,kW)&\mbox{by Lemma \ref{lemma:triangle}}\\
&&\\
&=\angle(kV,\mathcal{F})+\angle(V,W)&\\
&&\\
&=\angle(kV,\mathcal{F})+\angle(V,\mathcal{F})&\mbox{as\ }W=\pi_F(V)\\
&&\\
&\leq \angle(kV,\mathcal{F})+\angle(kV,\mathcal{F})&\mbox{by minimality}\\
&&\\
&=2\angle(kV,\mathcal{F})&
\end{array}
$$
and we are done.  
\end{proof}        

\subsection{Finishing the proof of the Eigenvalue Matching Theorem}

Armed with the lemmas of the previous two subsections, we now prove
Theorem \ref{prop:det-est}.

We begin by noting that the construction of $V'$ from $V$ above
respects subspace inclusion. I.e. if $U\subset V$ then $U'\subset V'$.
This follows from the definition of $V'$ and the fact that for two
singular subspaces $U_1$ and $U_2$ with $U_1\subset U_2$, we have
$K_{U_1}\cdot (W\cap K_{U_2}) \supset K_{U_2}$, where $W$ is the Weyl
group.

Now we simply proceed by induction on the number of vectors $k$. For
$k=1$ we set $V=v_1$ the statement of the proposition follows from
Lemma \ref{lem:subspace}.  Order the vectors by increasing angle with
$\mathcal{F}$. Assume the proposition for $k-1$ vectors, then set
$V_k=\op{span}\{v_1,\cdots,v_k\}$. By Lemma \ref{lem:subspace} we have
an orthogonal subspace of twice the dimension of $V_k$, namely $V_k'$,
which we may write by the preceeding paragraph as
$V'_k=V'_{k-1}\directsum W'$ where $W'$ is two dimensional. The same
lemma also guarantees that $\angle W',\mathcal{F}^\perp \leq C \angle
v_k,\mathcal{F}$, since $\angle v_k,\mathcal{F}=\angle
V_k,\mathcal{F}$.

This completes the proof of Theorem \ref{prop:det-est}.
         
\section{Finishing the proof of the Degree Theorem}
\label{section:finalproof}
We will break the proof of Theorem ~\ref{theorem:degree} into the
compact and noncompact cases.

\subsection{The compact case}

Suppose $M$ and $N$ are compact.  Since for $s>h(g)$, $F_s$ is a $C^1$
map, using Proposition \ref{prop:homotopy} and elementary integration
theory yields,

\begin{eqnarray}
  |\deg(f)| \vol(M)&=|\deg(f)| \int_{M} dg_{0}&=\left| \int_{N}f^*dg_{0}
\right| \cr  
&  =\left|\int_{N} F_s^*dg_{0} \right|&\cr 
& \leq\int_{N} \left| \Jac F_s\right| dg &\leq
C\left(\frac{s}{h(g_{0})}\right)^n \vol(N)
\end{eqnarray}

For the last inequality we have used the principal estimate from
Theorem~\ref{JacF}. Rearranging terms gives us the inequality in
Theorem~\ref{theorem:degree} since $C$ depends only on the dimension
and $\left(\frac{s}{h(g_{0})}\right)^n$ depends only on $n$ and the
smallest Ricci curvatures of $M$ and $N$.

\subsection{The noncompact case}
\label{sec:noncompact-case}

We now consider the case when $N$ (and/or $M$) has finite volume but is not
compact. In this setting, it is not known whether the limit in the
definition of $h(g)$ always exists. For this reason we will define the
quantity $h(g)$ to be
$$h(g)=\inf \left\{ s\geq 0 \left| \ \exists C>0 \text{ such that
      }\forall y\in Y, \ \int_Y e^{-s d(y,z)}
    dg(z)<C\right.\right\}.$$
In fact this agrees with the previous
definition for $h(g)$ when $N$ is compact. In the case of the
symmetric space $(M,g_0)$ this definition of $h(g_0)$ agrees with the
previous definition for compact manifolds.

For the finite volume case, the main difficulty is that, in order for
the proof given above to work, we need to know that $F_s$ is proper
(and thus surjective since $\op{deg}(F_s)=\op{deg}(f)\neq 0$).  For
this, we will need to prove higher rank analogs of some lemmas used in
\cite{BCS} for the rank one case. For the basics of degree theory for
proper maps between noncompact spaces, see \cite{FG}. We will need to
assume that the geometry of $N$ is bounded in the sense that its Ricci
curvatures are bounded from above and that the injectivity radius of
its universal cover $Y$ is bounded from below. These are the specific
assumptions implied in the third remark after the theorem.

We will show that $F_s$ is proper by essentially showing that the
barycenter of $\sigma^s_y$ lies nearby a convex set containing large
mass for this measure. This convex set is in turn far away from
$\phi(p)$ whenever $x$ is far from $p\in Y$.  We achieve this by first
estimating the concentration of the mass of $\sigma^s_y$ in certain
cones which will be our convex sets. One difficulty that arises in the
higher rank is that these cones must have a certain angle when
restricted to a flat. Another difficulty is that the ends of $M$ can
have large angle at infinity. In fact our methods breakdown unless we
control the asymptotic expansion of $f$ down the ends (see Remarks~
\ref{remarks:ends}).

First, we localize the barycenter of the measure $\sigma^s_y$. Let
$v_{(x,\theta)}$ be the unit vector in $S_xX$ pointing to $\theta\in\D
X$.

\begin{lemma}
\label{barycenter}
Let $K\subset X$ and $y\in Y$ be such that $(\phi_*\mu_y^s)(K)> C$ for
some constant $1>C>\frac12$. Suppose that for all $x\in X$ there
exists $v\in S_x X$ such that for all $z\in K$:
$$\int_{\D_F X}\inner{ v_{(x,\theta)},v}d\nu_z(\theta)\ge\frac1{C}-1$$ 
Then $$x\neq \tilde{F}_s(y)$$
\end{lemma}

\begin{proof}
If $\tilde{F}_s(y)=x$ then $\nabla_x \mathcal{B}_{s,y}(x)=0$.  However,
$\nabla_x \mathcal{B}_{s,y}(x)$ may be expressed as
$$\int_{X}\int_{\D_F
X}v_{(x,\theta)}d\nu_z(\theta)d\phi_*\mu_y^s(z)$$ 
where $v_{(x,\theta)}$
is the unit vector in $S_xX$ pointing to $\theta\in\D_F X$. Then we
have
\begin{align*}
  \left\Vert D_x\mathcal{B}_{s,y}\right\Vert
  &=\left\Vert\int_{X}\int_{\D_F
      X}  v_{(x,\theta)}d\nu_z(\theta)d\phi_*\mu_y^s(z) \right\Vert \\
  &\ge \left\Vert\int_K\int_{\D_F X}v_{(x,\theta)}d\nu_z
    (\theta)d\phi_*\mu_y^s(z)\right\Vert- \\
  &\hspace{2cm} \left\Vert\int_{X-K}\int_{\D_F X}
    v_{(x,\theta)}d\nu_z(\theta)d\phi_*\mu_y^s(z)\right\Vert \\
  &\ge \int_K\int_{\D_F X}\inner{v_{(x,\theta)},v}d\nu_z
  (\theta)d\phi_*\mu_y^s(z)-\phi_*\mu_y^s(X-K)\\
  &\ge \phi_*\mu_y^s(K)\left(\frac1{C}-1\right)-1+\phi_*\mu_y^s(K)\\
  &> C\left(\frac1{C}-1\right)-1+C=0\\
\end{align*}
The strictness of the inequality finishes the proof.
\end{proof}

For $v\in S X$ and $\alpha>0$ consider the convex cone,
$$E_{(v,\alpha)}=\exp_{\pi(v)}\left\{w\in T_{\pi(v)}X\ \vert\ 
  \angle_{\pi(v)} (v(\infty),w(\infty))\leq \alpha \right\},$$
where $\pi:TX \to X$ is the tangent bundle projection.

Denote by $\D E_{(v,\alpha)}\subset\D X$ its boundary at infinity.

\begin{lemma}
\label{onethird}
There exists $T_0>0$ and $\alpha_0>0$ such that for all $t\geq T_0$,
all $x\in X$, all $v\in S_xX$ and all $z\in E_{(g^{t}v,\alpha_0)}$,
$$\int_{\D_F X}\inner{v_{(x,\theta)},v}d\nu_z(\theta)\ge
\frac{\sqrt{2}}{3}.$$
\end{lemma}
\begin{proof}
  Since the isometry group of the symmetric space $X$ is transitive
  on $X$ and for any isometry $\psi$,
  $d\psi(E_{(v,\alpha)})=E_{(d\psi(v),\alpha)}$, it is sufficient to
  prove the lemma for a fixed $x$ and all $v\in S_xX$.
  
  For now choose $\alpha_0<\pi/4$. Take a monotone sequence $t_i\to
  \infty$, and any choice $z_i\in E_{(g^{t_i}v,\alpha)}$ for each
  $t_i$. It follows that some subsequence of the $z_i$, which we again
  denote by $\{z_i\}$, must tend to some point $\theta\in \D
  E_{(v,\alpha)}$.
  
  Let $\nu_\theta$ be the weak limit of the measures $\nu_{z_i}$.
  From Theorem~\ref{theorem:support}, $\nu_\theta$ is a probability
  measure supported on a set $S_\theta$ satisfying
$$\angle_x(\theta,\xi) \leq \frac{\pi}4 \quad \forall \xi\in
S_\theta.$$

Therefore we have,
\begin{align}
\label{eq:inner-product}
\int_{S_{\theta}}\inner{v_{(x,\xi)},v_{(x,\theta)}}d\nu_\theta(\xi)
\ge\frac{\sqrt{2}}{2}
\end{align}  

Now whenever $\theta\in \D E_{(v,\alpha)}$ then $v=v_{(x,\theta)}+\eps
v'$ for some unit vector $v'$ and $\eps\leq \sin(\alpha)$. Using
either case above we may write
$$\int_{\D_F X}\inner{v_{(x,\xi)},v}d\nu_\theta(\xi)\ge \int_{\D_F
  X}\inner{v_{(x,\xi)},v_{(x,\theta)}}d\nu_\theta(\xi)
-\sin(\alpha).$$
So choosing $\alpha$ small enough we can guarantee
that
\begin{enumerate}
\item any two Weyl chambers intersecting $E_{(g^tv,\alpha)}$ for all
  $t>0$ in the same flat must share a common face of dimension
  $\rank(M)-1$, and
\item for any $\theta\in \D E_{(v,\alpha)}$,
  $$\int_{\D_F X}\inner{v_{(x,\xi)},v}d\nu_\theta(\xi)\ge
  \frac{\sqrt{2}}{2.5}.$$
\end{enumerate}

Let
$$E_{(v(\infty),\alpha)}=\cap_{t>0}\D E_{(g^tv,\alpha)}.$$
By the first
property used in the choice of $\alpha$ above, for any two points
$\theta_1,\theta_2\in E_{(v(\infty),\alpha)}$, either $\theta_1$ and
$\theta_2$ are in the boundary of the same Weyl chamber, or else there
is another point $\theta'$ in the intersection of the boundaries at
infinity of the closures of the respective Weyl chambers. 

By maximality there is some $\theta_0\in E_{(v(\infty),\alpha)}$
intersecting the boundary at infinity of the closure of every Weyl
chamber which intersects $E_{(g^tv,\alpha)}$ for all $t>0$.  Hence,
for every $\theta \in E_{(v(\infty),\alpha)}$, the support of the
limit measure $\nu_\theta$ satisfies $S_\theta\subset S_{\theta_0}$.
(While $\theta_0$ is not necessarily unique, the support
$S_{\theta_0}$ of the corresponding limit measure $\nu_{\theta_0}$
is.)

As $t$ increases, for any $z\in E_{(g^tv,\alpha)}$, the measures $\nu_z$
uniformly become increasingly concentrated on $S_{\theta_0}$. Then
applying the estimate \eqref{eq:inner-product} to $\theta=\theta_0$,
we may choose $T_0$ sufficiently large so that for all $z\in
E_{(g^tv,\alpha)}$ with $t>T_0$,
$$\int_{\D_F X}\inner{v_{(x,\xi)},v}d\nu_z(\xi)\ge
\frac{\sqrt{2}}{3}.$$
\end{proof}

\begin{proposition}
$F_s$ is proper.
\end{proposition}

\begin{proof}
  By way of contradiction, let $y_i\in Y$ be an unbounded sequence
  such that $\{\tilde{F}_s(y_i)\}$ lies in a compact set $K$. We may
  pass to an unbounded subsequence of $\{y_i\}$, which we again denote
  as $\{y_i\}$, such that the sequence $\phi(y_i)$ converges
  within a fundamental domain for $\pi_1(M)$ in $X$ to a point
  $\theta_0\in \D X$. Since $K$ is compact, the set
  $$A=\bigcap_{x\in K}E_{(g^{T_0}v_{(x,\theta_0)},\alpha_0)}$$
  contains an open neighborhood of $\theta_0$ and $d_X(A,K)\geq T_0$.
  Notice that $A$ is itself a cone, being the intersection of cones on
  a nonempty subset of $\D X$.
  
  We now show that $A$ contains the image $\phi(B(y_i,R_i))$ of
  increasingly large balls ($R_i\to\infty$). However, we observe from
  the fact that $A$ is a cone on an open neighborhood of $\theta_0$ in
  $\D X$ that $A$ contains balls $B(\phi(y_i),r_i)$ with
  $r_i\to\infty$. By assumption $f$, and hence $\phi$, is coarsely
  Lipschitz:
  $$d_X(\phi x,\phi y)\leq K d_Y(x,y)+C$$
  for some constants $C>0$ and
  $K\geq 1$. Therefore $\phi^{-1}(B(\phi(y_i),r_i))\supset B(y_i,R_i)$
  where $K R_i+C> r_i$. In particular $R_i\to\infty$.
  
  Hence, there exists an unbounded sequence $R_i$ such that
  $B(y_i,R_i)\subset \phi^{-1}(A)$. Furthermore, since the Ricci
  curvature is assumed to be bounded from above and the injectivity
  radius from below, we have that $\vol (B(y_i,\op{injrad}))$ is
  greater than some constant independent of $y_i$ and hence $\int_Y
  e^{-s d(y_i,z)}dg(z)>Q$ for some constant $Q>0$. By choice of $s$
  there is a constant $C_s$ depending only on $s$ such that  $\int_Y
  e^{-s d(y,z)}dg(z)<C_s$ for all $y\in Y$.

  In polar coordinates we may write,
\begin{align*}
\int_Y e^{-s d(y,z)}dg(z)&=\int_0^\infty e^{-s t}\vol(S(y,t))dt\\
&=\int_0^\infty e^{-s t}\frac{d}{dt}\vol(B(y,t))dt\\
&=-\int_0^\infty \frac{d}{dt}\left(e^{-s t}\right)\ \vol(B(y,t))dt\\
&=s\int_0^\infty e^{-s t}\vol(B(y,t))dt.
\end{align*}
Using this we may estimate, using any $\delta<s-h(g)$,
\begin{align*}
  \mu_{y_i}^s(\phi^{-1}(A))&> \mu_{y_i}^s(B(y_i,R_i))\\
  &=1-\frac{\int_{R_i}^\infty e^{-s t}\vol(B(y_i,t))dt}{\int_0^\infty
    e^{-s
      t}\vol(B(y_i,t))dt}\\
  &\geq 1-\frac{e^{-\delta R_i}\int_{R_i}^\infty e^{-(s-\delta)
      t}\vol(B(y_i,t))dt}{\int_0^\infty e^{-s
      t}\vol(B(y_i,t))dt}\\
  &\geq 1-e^{-\delta R_i}\frac{C_{s-\delta}}{Q}.
\end{align*}

Therefore for all sufficiently large $i$,
$$\mu_{y_i}^s(\phi^{-1}(A))> \frac{3}{3+\sqrt{2}}.$$
The constant
$\frac{3}{3+\sqrt{2}}$ is the constant $C$ from Lemma \ref{barycenter}
such that $\frac{1}{C}-1=\frac{\sqrt{2}}{3}$.
  
  Set $v_i=g^{T_0+1}v_{(\tilde{F}_s(y_i),\theta_0)}$. Recalling that
  $A\subset E_{(v_i, \alpha_0)}$ for all $i$, we have that for
  sufficiently large $i$,
  $$\phi_*\mu_{y_i}^s(E_{(v_i, \alpha_0)})> \frac{3}{3+\sqrt{2}}$$
  but
  $d_X(\tilde{F}_s(y_i),E_{(v_i,\alpha_0)})>T_0$, contradicting the
  conclusion of Lemma \ref{barycenter} in light of Lemma
  \ref{onethird}.
\end{proof}

\begin{remarks}
\label{remarks:ends}
\indent
\begin{enumerate}
\item In the proof of the above proposition, we used that
  $\op{injrad}$ is bounded from below and $\op{Ricci}$ curvature is
  bounded from above only to show that the volume of balls of any
  fixed radius are bounded from below.
  
\item Ideas from coarse topology can be used to remove the coarse
  lipschitz assumption on $f$ in the case that the ends of $M$ have
  angle at infinity bounded away from $\pi/2$. However, $M$ may
  have ends containing pieces of flats with wide angle (consider the
  product of two rank one manifolds each with multiple cusps, or for a
  complete classification of higher rank ends see \cite{Ha}). For
  such spaces it is possible to construct a proper map $f:M\to M$ such
  that for a radial sequence $y_i\to \infty$, $\phi$ maps the bulk of
  the mass of $\mu^s_{y_i}$ into a set (almost) symmetrically arranged
  about the point $p\in X$ thus keeping $\tilde{F}_s(y_i)$ bounded.
  This explains the need for a condition on $f$ akin to the coarse
  lipschitz hypothesis.
\end{enumerate}
\end{remarks}

The inequality in Theorem~\ref{theorem:degree} now follows as in the
compact case, with $\deg(f)$ and $\deg(F_s)$ suitably interpreted.

\providecommand{\bysame}{\leavevmode\hbox to3em{\hrulefill}\thinspace}

\noindent
Christopher Connell:\\
Dept. of Mathematics, University of Illinois at Chicago\\
Chicago, IL 60680\\
E-mail: cconnell@math.uic.edu
\medskip

\noindent
Benson Farb:\\
Dept. of Mathematics, University of Chicago\\
5734 University Ave.\\
Chicago, Il 60637\\
E-mail: farb@math.uchicago.edu

\end{document}